\documentclass[11pt,twoside]{article}
\usepackage{latexsym}
\usepackage{amssymb,amsbsy,amsmath,amsfonts,amssymb,amscd}
\usepackage{epsfig, graphicx}
\usepackage{graphicx}
\usepackage{subfigure}
\usepackage{color}
\usepackage{epstopdf}
\newcommand{\commentout}[1]{}

\newcommand {\Chi} {{\bf \raise 2pt \hbox{$\chi$}} }
\newcommand {\f}   {\frac}
\newcommand {\p}   {\partial}

\newcommand {\proof} {\noindent {\bf Proof}. }
\newcommand{\beq}{\begin{equation}}
\newcommand{\eeq}{\end{equation}}
\newcommand{\bea} {\begin{array}{rl}}
\newcommand{\eea} {\end{array}}
\newcommand{\bepa}{\left\{ \begin{array}{l}}
\newcommand{\eepa} {\end{array}\right.}
\newtheorem{theorem}{Theorem}[section]

\newtheorem{remark}[theorem]{Remark}

\newcommand{\qed}{{ \hfill
                       {\unskip\kern 6pt\penalty 500 \raise -2pt\hbox{\vrule\vbox to 6pt{\hrule width 6pt
                       \vfill\hrule}\vrule} \par}   }}
\title{\Large \bf An Asymptotic Preserving method for strongly anisotropic diffusion equations based on field line integration}

\author{Min Tang
\thanks{Institute of natural sciences and department of mathematics,
Shanghai Jiao Tong University, Shanghai, 200240, China.}
\and Yihong Wang
\thanks{Institute of natural sciences and department of mathematics,
Shanghai Jiao Tong University, Shanghai, 200240, China.}}

\date{\today}

\begin{document}
\maketitle
\pagestyle{plain}
\pagenumbering{arabic}

\begin{abstract}
In magnetized plasma, the magnetic field confines the particles around the field lines. The anisotropy intensity in the viscosity and heat conduction may reach the order of $10^{12}$.
When the boundary conditions are periodic or Neumann, the strong diffusion leads to an ill-posed limiting problem.
To remove the ill-conditionedness in the highly anisotropic diffusion equations, we introduce a simple but very efficient asymptotic preserving
reformulation in this paper.
The key idea is that, instead of discretizing the Neumann boundary conditions locally, we
replace one of the Neumann boundary condition by the integration of the original problem along the
 field line, the singular $1/\epsilon$ terms can be replaced by $O(1)$ terms after the integration, so that
 yields a well-posed problem.  Small modifications to the original code are required and no change of coordinates nor mesh adaptation are needed.
 Uniform convergence with respect to the anisotropy strength $1/\epsilon$ can be observed numerically and
the condition number does not scale with the anisotropy.

 \end{abstract}
\\[3mm]
 \section{Introduction}\label{intro}
Anisotropic diffusion is encountered in many physical applications, including flows in porous media \cite{Ashby 99},
heat conduction in fusion plasmas \cite{Bram 14},
atmospheric or oceanic flows \cite{Galperin 10,Treguier 02} and so on.
In magnetized plasmas, since the particles
are confined by the magnetic field, along the field lines, the distance between two successive collision is much larger than
the mean free path in the perpendicular direction, which yields extremely anisotropic diffusion tensors and
their directions may vary in space and time. Numerically, the error in the diffusion
parallel to the magnetic field may pollute the small perpendicular diffusion. The numerical resolution of highly anisotropic physical problems is a challenging task.
As summarized in \cite{Bram 14}, several problems may arise:
1) yield anisotropy
dependent convergence \cite{Babuska 92}; 2) introduce large perpendicular errors \cite{Umansky 05};
3) loss positivity near high gradients.

There is another difficulty which relates to specific boundary conditions.
As discussed in \cite{Degond122}, when the boundary conditions are periodic or Neumann, the strong diffusion leads to an ill-posed limiting problem.
This difficulty arises only for Neumann or periodic boundary conditions as well as magnetic re-connections, while vanishes for Dirichlet and Robin boundary conditions.
However,  a lot of physical problems, for instance the tokamak plasmas on a torus, atmospheric plasma, periodic or Neumann boundary conditions have a considerable impact \cite{Degond122}.


The problem under consideration is the following two-dimensional advection diffusion equation with anisotropic diffusivity:
\beq\label{eq:ellipT}
\left\{\begin{array}{ll}
-\nabla\cdot (A\nabla u^{\epsilon})=f,&\mbox{in $\Omega$},\\
 \mathbf{n}\cdot A\nabla u^{\epsilon}=0 ,&\mbox{on $\Gamma_N$},\\
u^{\epsilon}=0 ,&\mbox{on $\Gamma_D$},
\end{array}
\right.
\eeq
where $\Omega\subset R^2$ is a bounded domain and $\mathbf{n}$ is the outward normal vector.
The boundary $\Gamma=\partial \Omega$ is composed by two parts $\Gamma_N\cup \Gamma_D$, one with Neumann boundary condition and the other with Dirichlet.
The direction of the
anisotropy or the magnetic field line is given by a unit vector field $\mathbf{b}=(\cos\theta,\sin\theta)$.
 The two parts of the boundary relate to the magnetic vector field and are given by $\Gamma_D:=\{x\in\partial\Omega\mid \mathbf{b}(x)\cdot \mathbf{n}=0\}$
 and $\Gamma_N:=\partial \Omega\setminus \Gamma_D$. The anisotropic
diffusion matrix is given by
 \beq\label{eq:A}
 A(x,y)=\left(\begin{array}{cc}\cos\theta&-\sin\theta\\ \sin\theta&\cos\theta\end{array}\right)
\left(\begin{array}{cc}1/\epsilon&0\\0&\alpha\end{array}\right)\left(\begin{array}{cc}\cos\theta&\sin\theta\\-\sin\theta&\cos\theta\end{array}\right),
\eeq where $\alpha$ is at $O(1)$ and the parameter $0<\epsilon<1$ can be very small. All coefficients $\alpha$, $\epsilon$ and $\theta$ can depend on time and space.
The problem becomes high anisotropy when $\epsilon\ll1$. Let $(\xi,\eta)$ be the aligned coordinate system,
the formal limit of $\epsilon\to 0$ leads to
$$-\p^2_\xi u=0,\qquad \p_\xi u=0,\quad \mbox{on $\Gamma_N$},\qquad u=0,\quad\mbox{on $\Gamma_D$}.
$$
Any function constant along the $\mathbf{b}$ field solve the above equation, thus there exist infinitely many solutions.
 Consequently, the condition number of the discretized system
 tends to $\infty$ and the solution will suffer from round-off errors.

This difficulty has been observed and studied in a series paper by Degond et.al  \cite{Degond101,Degond102,Degond121, Degond122,Negulescu14}.
Their methods are refereed as asymptotic preserving (AP) method in the sense that they can deal with all cases that $\epsilon$ ranges from $O(1)$ to very small.
AP methods can capture the limiting solution using coarse meshes, they have been successfully developed for various models, we refer to the lecture notes by Jin in \cite{Jin 10}.
For the strongly anisotropic diffusion equation considered in this paper, the AP methods have the advantage that the condition number does not scale with the anisotropy, thus no preconditioner is required.
The main idea is to decompose the solution into two parts: a mean part along the anisotropic direction and a fluctuation part.
Then reformulate the original equation into a modified system which, as $\epsilon\to 0$, will reduce to a well-posed system satisfied by the limiting solution.
The special case when the anisotropy direction is aligned with one of the coordinate is studied in  \cite{Degond101}. Extensions to arbitrary fields are proposed in \cite{Degond102},
where Cartesian grids are used and no change of coordinate, no mesh adaptation are required. However the method in \cite{Degond102} results in a considerably bigger linear system. In order
to reduce the computational cost, a modification is presented in \cite{Degond122}. All these methods are based on solution splitting and equation reformulation.

In this paper, we propose a different approach that is simple and easily extendable. The key idea is that:
The condition $\mathbf{b}\cdot\mathbf{n}\neq 0$ indicates that the boundary is not parallel to
the field line, thus at all points belong to $\Gamma_N$, there exists one field line cross the point. The two boundary conditions for each field line are Neumann
according to the definition of  $\Gamma_N$ and $\Gamma_D$. Therefore, for each point $(x_1,y_1)$ belongs \textcolor{blue}{to} $\Gamma_N$, the field line cross this point will connect to another point $(x_2,y_2)$
also belongs to $\Gamma_N$. Neumann boundary condition $\mathbf{n}\cdot A\nabla u^{\epsilon}=0$ is satisfied at both $(x_1,y_1)$ and $(x_2,y_2)$. The fact that the limiting problem becomes
ill-posed can be understood as the existence of infinitely many solutions to the one dimensional problem $-\p^2_\xi u=0$ with $\p_\xi u=0$ at the boundary ($\xi$ is the local coordinate along the field line). Therefore,
 instead of discretizing the two Neumann boundary conditions locally, we
replace one of the Neumann boundary condition by the integration of the original problem along the
 field line, the singular $1/\epsilon$ terms can be replaced by $O(1)$ terms after the integration, so that
 yields a well-posed problem. The novelties of our scheme are listed below:
\begin{itemize}
\item Small modifications to the original code are required, which makes it attractable to engineers.
The idea can be coupled with most standard discretizations and the computational cost keeps almost the same.

\item Extensions to space dependent $\alpha$, $\epsilon$ and $\theta$ in \eqref{eq:A} are straight forward.
Uniform convergence with respect to the anisotropy strength $\epsilon$ can be observed numerically and the condition number does not scale with the anisotropy.
 \item No change of coordinates nor mesh adaptation are required.
 In the context of tokamak plasma simulation, the anisotropy being driven by the magnetic field which is time dependent.
The computational cost of magnetic field aligned coordinate is expensive which motivates the use of coordinates and meshes that are independent of the anisotropic direction \cite{Gunter05}.
\end{itemize}

Numerical simulations for anisotropic diffusivity attract a lot of researchers and engineers, see the review in \cite{Herbin 08}. For strong anisotropic problem, since the numerical error in the diffusion
parallel to the magnetic field may pollute the small perpendicular diffusion.
Common approach is to use magnetic field aligned coordinates, which may run into problems when there are magnetic re-connections or
highly fluctuating field directions.
Schemes used today include Finite volume method \cite{Potier 05,Sheng 09, Yuan 08}, Finite difference method \cite{Bram 14}, Mimetic Finite difference method \cite{Hyman 02}, Discontinuous Galerkin method \cite{ARNOLD82,ALEXANDRE 09},
Finite element method \cite{Gunter07, Hou 97, Li 10, Pasdunkorale 05} and so on. In \cite{Umansky 05}, the authors discuss about the numerical pollution in strong anisotropic problem, which is hard to avoid
for schemes using Cartesian grids. However pollution is
not the goal of our present paper, the standard discretization is enough to illustrate the idea of how to reformulate the ill-posed problem into a well-posed one.
Therefore, we use the simplest $5$-point (when the coordinate is aligned with the field line) or $9$-point
(when the coordinate is misaligned with the field line) finite difference method.

The paper is organized as follows. In section 2, we illustrate the idea of the reformulation in the special case that the magnetic field line is aligned with one of
the coordinates, then extend it to nonaligned case. The numerical discretizations are given in Section 3, where
 five and nine point finite difference method (FDM) are used respectively for the aligned and nonaligned case.
 Several numerical examples are presented in section 4 and the uniform convergence as well as the upper bound of the condition number, with respect to the anisotropy, can be observed. Finally we conclude with some discussions in section 5.

 \section{ The Asymptotic Preserving Reformulation}\label{num}

Let $\mathbf{n}$ be the unit outward normal on the boundary $\Gamma$. We decompose $\Gamma$ into three parts by the sign of $\mathbf{n}\cdot\mathbf{b}(x)$ such that:
$$\Gamma_D:=\{x\in\Gamma\mid \mathbf{n}\cdot\mathbf{b}(x)=0\},\quad \Gamma_{in}:=\{x\in\Gamma\mid \mathbf{n}\cdot\mathbf{b}(x)<0\},
\quad \Gamma_{out}:=\{x\in\Gamma\mid \mathbf{n}\cdot\mathbf{b}(x)>0\}.$$
It is obvious that $\Gamma_{N}:=\Gamma_{in}\cup\Gamma_{out}$.
\subsection{The special case when $\theta\equiv 0$}
 When $\theta\equiv 0$ in the diffusion matrix $A(x,y)$,  \eqref{eq:ellipT}  becomes
\beq\label{eq:ellipT1}
\left\{\begin{array}{ll}
\displaystyle - \partial_x\big(\frac{1}{\epsilon}\partial_xu^{\epsilon}\big)
-\partial_y\big(\alpha\partial_ yu^\epsilon\big)=f(x,y),\quad \mbox{in $\Omega$},\\
\displaystyle \partial_xu^{\epsilon}=0,\quad\mbox{on $\Gamma_N$},\qquad u^{\epsilon}=0, \quad\mbox{on $\Gamma_D$}.
\end{array}
\right.
\eeq
The formal limit of $\epsilon\rightarrow 0$ in  \eqref{eq:ellipT1} yields a ill-posed problem
\beq\label{eq:ellipT2}
-\partial_x^2 u^{\epsilon}=0,\quad \mbox{in $\Omega$},\qquad
 \partial_x u^{\epsilon}=0,\quad\mbox{on $\Gamma_N$},\qquad
 u^{\epsilon}=0,\quad\mbox{on $\Gamma_D$}.
\eeq
We consider a rectangular domain $\Omega=[0,a]\times[0,b]$ and that the field line is parallel to the $x$ axis. The Neumann boundary conditions are imposed on the left and right boundary.
\eqref{eq:ellipT1} is equivalent to the following system
\beq\label{eq:ellipT3}
\left\{\begin{array}{ll}
\displaystyle - \partial_x\big(\frac{1}{\epsilon}\partial_xu^{\epsilon}\big)
-\partial_y\big(\alpha\partial_ yu^\epsilon\big)=f(x,y), &\mbox{in $[0,a]\times[0,b]$},\\
\displaystyle\partial_x u^{\epsilon}=0,&\mbox{on $0\times[0,b]$},\\
\displaystyle  -\int_0^a\partial_y\big(\alpha\partial_ yu^\epsilon\big)\,dx=\int_0^af(x,y)\,dx,&\mbox{on $a\times[0,b]$},\\
\displaystyle  u^{\epsilon}=0, &\mbox{on $[0,a]\times 0\cup [0,a]\times b$}.
\end{array}
\right.
\eeq
The limiting problem of the above equation, when $\epsilon\to 0$, becomes
\beq\label{eq:ellipT3l}
\left\{\begin{array}{ll}
\displaystyle \p^2_xu^{0} =0,&\mbox{in $\Omega$},\\
\displaystyle\p_xu^{0}=0,&\mbox{on $0\times[0,b]$},\\
\displaystyle  -\int_0^a\partial_y\big(\alpha\partial_ yu^0\big)\,dx=\int_0^af(x,y)\,dx,&\mbox{on $a\times[0,b]$},\\
\displaystyle  u^{0}=0, &\mbox{on $[0,a]\times 0\cup [0,a]\times b$}.
\end{array}
\right.
\eeq
The solution to \eqref{eq:ellipT3l} can be easily obtained. First of all, the first two equations yield that, for fixed $y$, $u^0(x,y)$ is a constant for all $x\in[0,a]$.
Then let $ \int_0^au^0\,dx=\bar{u}^0(y)$, the third and forth equations in \eqref{eq:ellipT3l} give
\beq\label{eq:limitaligned}
 -\int_0^a\partial_y\big(\alpha\partial_ yu^0\big)\,dx=-\p_y\Big(\big(\f{1}{a}\int_0^a\alpha\,dx\big)\p_y\bar{u}^0\Big)=\int_0^af(x,y)\,dx,\qquad \bar{u}^0(0)=\bar{u}^0(b)=0,
\eeq
which is exactly the same as the limiting model in \cite{Degond101}. \eqref{eq:limitaligned} determines a unique solution $\bar{u}^0(y)$, together with the fact that
$u^0(x,y)$ is independent of $x$, \eqref{eq:ellipT3l} is consequently a well-posed problem.

 The format of this simple case is close to the scheme in \cite{Degond101}, where the AP property is achieved by macro-micro decomposition,
the average of $u^\epsilon$ satisfies a similar equation as the third equation in \eqref{eq:ellipT3}.
However, the conceptual idea of scheme design is entirely different.
We would like to
emphasis that, the case of $\theta=0$ is used to illustrate the basic idea of "field line integration". We can easily see that the ill-posed limiting problem becomes well-posed after replacing one of the boundary conditions by field line integration. Our idea of "field line integration" can be easily extended to more general cases.

\subsection{The general case}
 Let $\mathbf{b}=(\cos\theta,\sin\theta)$, $\mathbf{b}_\bot=(-\sin\theta,\cos\theta)$,  \eqref{eq:ellipT} can be written as
\beq\label{eq:SimellipA}
\left\{\begin{array}{ll}
\displaystyle  -  ( \mathbf{b} \cdot \nabla) (\frac{1}{\varepsilon}  \mathbf{b}\cdot \nabla u^{\epsilon}) -
     ( \nabla \cdot \mathbf{b}^T) ( \frac{1}{\epsilon}\mathbf{b} \cdot \nabla u^{\epsilon}) - \nabla
   \cdot \big( \alpha \mathbf{b}_{\perp} (\mathbf{b}_{\perp}\cdot \nabla u^{\epsilon})\big) = f ,& \mbox{in $\Omega$},\\
\displaystyle   \frac{1}{\epsilon} \mathbf{n}\cdot\mathbf{b}^T  (\mathbf{b}\cdot \nabla u^{\epsilon}) + \alpha \mathbf{n}\cdot\mathbf{b}^T_{\perp}
   (\mathbf{b}_{\perp}\cdot \nabla u^{\epsilon}) = 0 , &\mbox{on $\Gamma_N$},\\
\displaystyle u^\epsilon=0,&\mbox{on $\Gamma_D$}.
\end{array}
\right.
\eeq

We can pick any $\mathbf{b}$-field line $l$ and parameterize it by the arc length $s$, so that
$s = 0$ corresponds to a point on $\Gamma_{{in}}$ and $s = L_l$ (with
$L_l$ being the length of $l$) corresponds to a point on
$\Gamma_{{out}}$. Accordingly, $\partial_s$ will denote the derivative
along the line $l$, i.e. $\partial_s = \mathbf{b} \cdot \nabla$. Then, the first equation in \eqref{eq:SimellipA}
can be written on $l$ as

\beq\label{eq:rew1}
 -  \partial_s(\frac{1}{\epsilon}\p_s u^{\epsilon}) -  ( \nabla
   \cdot \mathbf{b}^T) \frac{1}{\epsilon} \partial_s^{} u^{\epsilon} -  \nabla
   \cdot \big( \alpha \mathbf{b}_{\perp} (\mathbf{b}_{\perp}\cdot \nabla u^{\epsilon})\big)  = f.
   \eeq
   In the strong anisotropic diffusion limit such that $\epsilon\rightarrow 0$, \eqref{eq:SimellipA} yields
\beq\label{eq:rew2lim}
- \partial_s^2 u^0 -  ( \nabla
   \cdot \mathbf{b}^T) \partial_s^{} u^0 =0,\quad \mbox{in $\Omega$},\qquad
 \partial_s^{} u^{0}=0,\quad\mbox{on $\Gamma_N$},\qquad
 u^{0}=0,\quad\mbox{on $\Gamma_D$}.
\eeq
 Consider the field lines whose starting and ending points belong to $\Gamma_N$, as any function that keeps constant along those field lines satisfies \eqref{eq:rew2lim}, there is no uniqueness of the solution to the limiting model \eqref{eq:rew2lim}. However, similar as in
 \eqref{eq:ellipT3}, we
 can reformulate \eqref{eq:SimellipA} by replacing the Neumann boundary condition on one side of each field line by the integration of the equation,
 so that to obtain a well-posed limiting problem. The details are as follows.

First of all, we introduce a function $E$ which solves on $l$ the differential equation
\beq\label{eq:E}  \partial_s E = ( \nabla \cdot \mathbf{b}^T) E    \eeq
with initial condition $E = 1$ at $s = 0$.
Then we can combine the two singular $O(1/\varepsilon)$ terms in \eqref{eq:rew1} and rewrite them into a conservative form such that
\beq\label{eq:rew2}
-  \partial_s^{} ( E \frac{1}{\epsilon}\partial_s^{} u^{\epsilon}) - E \nabla
   \cdot \big( \alpha \mathbf{b}_{\perp} (\mathbf{b}_{\perp}\cdot \nabla u^{\epsilon})\big)  = {Ef}.
   \eeq
Integrating over $l$ and using the boundary conditions on $\Gamma_{{in}}$
and $\Gamma_{{out}}$, which for $s = 0$ and $s = L_l$ take the form:
\[ \frac{1}{\epsilon} \mathbf{n}\cdot\mathbf{b}^T  \p_s u^{\epsilon} + \alpha \mathbf{n}\cdot\mathbf{b}^T_{\perp}
   (\mathbf{b}_{\perp}\cdot \nabla u^{\epsilon}) = 0 \]
gives

\beq\label{eq:rew3}
\left[ E \alpha \frac{\mathbf{n} \cdot \mathbf{b}^T_{\perp}}{\mathbf{n} \cdot \mathbf{b}^T} \mathbf{b}_{\perp}^{} \cdot
   \nabla u^{\epsilon} \right]_{s = 0}^{s = L_l} - \int_0^{L_l} E \nabla
   \cdot \big( \alpha \mathbf{b}_{\perp} (\mathbf{b}_{\perp}\cdot \nabla u^{\epsilon})\big)  {ds} = \int_0^{L_l} Ef{ds}.
\eeq

Therefore, \eqref{eq:SimellipA} can be rewritten
 as
  \beq\label{eq:SimellipC}
\left\{\begin{array}{ll}
\displaystyle  -  ( \mathbf{b} \cdot \nabla) (\frac{1}{\varepsilon}  \mathbf{b}\cdot \nabla u^{\epsilon}) -
     ( \nabla \cdot \mathbf{b}^T) ( \frac{1}{\epsilon}\mathbf{b} \cdot \nabla u^{\epsilon}) - \nabla
   \cdot \big( \alpha \mathbf{b}_{\perp} (\mathbf{b}_{\perp}\cdot \nabla u^{\epsilon})\big) = f ,& \mbox{in $\Omega$},\\
\displaystyle   \frac{1}{\epsilon} \mathbf{n}\cdot\mathbf{b}^T  (\mathbf{b}\cdot \nabla u^{\epsilon}) + \alpha \mathbf{n}\cdot\mathbf{b}^T_{\perp}
   (\mathbf{b}_{\perp}\cdot \nabla u^{\epsilon}) = 0 , &\mbox{on $\Gamma_{in}$},\\
\displaystyle \left[ E \alpha \frac{\mathbf{n} \cdot \mathbf{b}^T_{\perp}}{\mathbf{n} \cdot \mathbf{b}^T} \mathbf{b}_{\perp}^{} \cdot
   \nabla u^{\epsilon} \right]_{s = 0}^{s = L_l} - \int_0^{L_l} E \nabla
   \cdot \big( \alpha \mathbf{b}_{\perp} (\mathbf{b}_{\perp}\cdot \nabla u^{\epsilon})\big)  {ds} = \int_0^{L_l} Ef{ds}.
&\mbox{on $\Gamma_{out}$},
\\
u^\epsilon=0. &\mbox{on $\Gamma_D$}.
\end{array}
\right.
\eeq
The above system is equivalent to the original system \eqref{eq:ellipT}, so that is well-posed for fixed $\epsilon\neq 0$.
When $\epsilon\to 0$, since $\mathbf{n}\cdot\mathbf{b}^T\neq 0$ for $\forall (x,y)\in\Gamma_{out}\subset\Gamma_N$, the leading order of \eqref{eq:SimellipC} gives
  \beq\label{eq:SimellipD}
\left\{\begin{array}{ll}
\displaystyle  - \partial_s^2 u^0 - ( \nabla
   \cdot \mathbf{b}^T) \partial_s^{} u^0 = 0,& \mbox{in $\Omega$},\\
\displaystyle  \partial_s u^0 = 0 ,  &\mbox{on $\Gamma_{in}$},\\
\displaystyle \left[ E \alpha \frac{\mathbf{n} \cdot \mathbf{b}^T_{\perp}}{\mathbf{n} \cdot \mathbf{b}^T} \mathbf{b}_{\perp}^{} \cdot
   \nabla u^{0} \right]_{s = 0}^{s = L_l} - \int_0^{L_l} E \nabla
   \cdot \big( \alpha \mathbf{b}_{\perp} (\mathbf{b}_{\perp}\cdot \nabla u^{0})\big)  {ds} = \int_0^{L_l} Ef{ds}.
&\mbox{on $\Gamma_{out}$},
\\
u^0=0. &\mbox{on $\Gamma_D$}.
\end{array}
\right.
\eeq
The first two equations in \eqref{eq:SimellipD} give that $u^0$ is a constant along $l$ while the third equation yields
the equation to determine the constant.

The reformulated equation \eqref{eq:rew2} is only for derivation of the third equation in \eqref{eq:SimellipC}.
The numerical discretizations can be performed on the Cartesian coordinates $(x,y)$.
Any standard algorithms can be applied to
discretize the first equation in \eqref{eq:SimellipC}, i.e. the first equation in \eqref{eq:ellipT} inside $\Omega$.
 The integration along $l$ with respect to $s$ in the third equation
of \eqref{eq:SimellipC} can be performed as follows:

First of all, we calculate $E$ that satisfies \eqref{eq:E}, whose solution can be given explicitly:
 \beq E=e^{\int_0^ s \nabla \cdot \mathbf{b}^T  \,\,ds}.\label{eq:Exieta}\eeq
Since $\mathbf{b}$ is given, so is $\nabla\cdot\mathbf{b}$, we can use  trapezoidal's rule to discretize the integration along the field line.
 Then,
from
\beq
 - \nabla
   \cdot \big( \alpha \mathbf{b}_{\perp} (\mathbf{b}_{\perp}\cdot \nabla u^{\epsilon})\big) =-\nabla\cdot\Big(\left(\begin{array}{cc}\cos\theta&-\sin\theta\\ \sin\theta&\cos\theta\end{array}\right)
\left(\begin{array}{cc}0&0\\0&\alpha\end{array}\right)\left(\begin{array}{cc}\cos\theta&\sin\theta\\-\sin\theta&\cos\theta\end{array}\right)\nabla u^{\epsilon}\Big),
\label{eq:alphamatrix}\eeq
we only need to
find the approximation to
$$Ef+E\nabla\cdot\Big(\left(\begin{array}{cc}\cos\theta&-\sin\theta\\ \sin\theta&\cos\theta\end{array}\right)
\left(\begin{array}{cc}0&0\\0&\alpha\end{array}\right)\left(\begin{array}{cc}\cos\theta&\sin\theta\\-\sin\theta&\cos\theta\end{array}\right)\nabla u^{\epsilon}\Big)
,
$$
on each integration point and then use trapezoidal's rule to approximate the integral in the third equation
of \eqref{eq:SimellipC}. The obtained integration is equal to a discretization of $\left[ E \alpha \frac{\mathbf{n} \cdot \mathbf{b}^T_{\perp}}{\mathbf{n} \cdot \mathbf{b}^T} \mathbf{b}_{\perp}^{} \cdot
   \nabla u^{0} \right]_{s = 0}^{s = L_l}$. More details about the discretization will be discussed in section 3.

\section{The Numerical Discretization}
For simplicity, we consider a rectangular domain and uniform Cartesian grids in this present paper. The notations for the domain
and grid nodes are
$\Omega=[0,a]\times[0,b]$ and
 $$
\begin{array}{ll}
\mathbf{z}_{i,j}=(x_i,y_j),\qquad&\mbox{with $i=0,1,\cdots I$; $j=0,1,\cdots J$}.
\end{array}
$$
Here $I$ and $J$ are two positive integers,
$x_i=ih_x$, $y_j=jh_y$ with $h_x=a/I$, $h_y=b/J$.
{To simplify the notations, we consider
the left and right boundaries belong to $\Gamma_N$ and
the bottom and top boundaries belong to $\Gamma_D$.}

\subsection{$5$-point FDM for the aligned case}
\begin{figure}[htb]
\centering
{
\includegraphics[width=7.5cm] {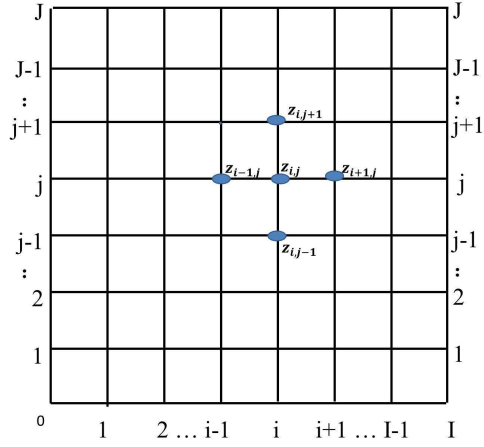}}
\caption{
 five point stencil.}\label{fig:5point}
\end{figure}

Let us begin with the aligned case \eqref{eq:ellipT1}.
For the $j$th row along the $y$ direction in Figure  \ref{fig:5point},
we can
use the standard five point FDM at the internal points $(x_i,y_j)$ such that
\beq\label{eq:ellipsT1scheme1}
\begin{aligned}
&-\f{1}{h_x}\Big(\f{1}{2}\big(\f{1}{\epsilon_{i+1,j}}+\f{1}{\epsilon_{i,j}}\big)\f{u^\epsilon_{i+1,j}-u^\epsilon_{i,j}}{h_x}-
\f{1}{2}\big(\f{1}{\epsilon_{i,j}}+\f{1}{\epsilon_{i-1,j}}\big)\f{u^\epsilon_{i,j}-u^\epsilon_{i-1,j}}{h_x}\Big)
\\&\qquad-
\f{1}{h_y}\Big(\f{1}{2}\big(\alpha_{i,j+1}+\alpha_{i,j}\big)\f{u^\epsilon_{i,j+1}-u^\epsilon_{i,j}}{h_y}-
\f{1}{2}\big(\alpha_{i,j}+\alpha_{i,j-1}\big)\f{u^\epsilon_{i,j}-u^\epsilon_{i,j-1}}{h_y}\Big)=f_{i,j}.
\end{aligned}
\eeq
In the limit $\epsilon\to 0$,
\eqref{eq:ellipsT1scheme1} yields
\beq\label{eq:ellipsT1scheme2}
u^0_{i+1,j}-u^0_{i,j}=u^0_{i,j}-u^0_{i-1,j},\qquad \mbox{for $i=1,\cdots,I-1$}.
\eeq
If we discretize the Neumann boundary condition on the left boundary by
$u^\epsilon_{0,j}-u^\epsilon_{1,j}=0$, together with \eqref{eq:ellipsT1scheme2}, we find $u^0_{I,j}-u^0_{I-1,j}=0$. Therefore,
if we use $u^\epsilon_{I-1,j}=u^\epsilon_{I,j}$ to discretize the Neumann boundary condition on the right, the coefficient matrix becomes ill-conditioned.

Due to the above observation, in order to remove the ill-posedness, we have to discretize the Neumann boundary condition on the right boundary in a different way.
The third equation in \eqref{eq:ellipT3} is used for each $y_j$ and we discretized it by the simplest first order approximation such that
\beq\label{eq:ellipsT1scheme4}
-\sum_{i=1}^{I-1}
\f{1}{h_y}\Big(\f{1}{2}\big(\alpha_{i,j+1}+\alpha_{i,j}\big)\f{u^\epsilon_{i,j+1}-u^\epsilon_{i,j}}{h_y}-
\f{1}{2}\big(\alpha_{i,j}+\alpha_{i,j-1}\big)\f{u^\epsilon_{i,j}-u^\epsilon_{i,j-1}}{h_y}\Big)
=\sum_{i=1}^{I-1} f_{i,j}.
\eeq
 The following theorem shows that the resulting discretized system become well-posed for any $\epsilon>0$.
\begin{theorem}\label{theo:unique}
The discretized scheme that using $5$-point FDM in the internal grids and $u^\epsilon_{0,j}-u^\epsilon_{1,j}=0$ on the left boundary, \eqref{eq:ellipsT1scheme4} on the right boundary
is uniquely solvable for any $\epsilon\geq 0$.
\end{theorem}
\proof
Let
$$\begin{aligned}
\kappa_{i,j}^r=\f{1}{2h_x^2}\big(\f{1}{\epsilon_{i+1,j}}+\f{1}{\epsilon_{i,j}}\big),\quad& \kappa_{i,j}^l=\f{1}{2h_x^2}\big(\f{1}{\epsilon_{i-1,j}}+\f{1}{\epsilon_{i,j}}\big),\\
\kappa_{i,j}^t=\f{1}{2h_y^2}\big(\alpha_{i,j+1}+\alpha_{i,j}\big),\quad& \kappa_{i,j}^b=\f{1}{2h_y^2}\big(\alpha_{i,j-1}+\alpha_{i,j}\big).
\end{aligned}
$$
After rearrangement, for $\forall i\in\{1,\cdots, I-1\}$, the $5$-point FDM \eqref{eq:ellipsT1scheme1} can be written as
\beq\label{eq:5dfminternal}
\kappa_{i,j}^tu^\epsilon_{i,j+1}
-\big(\kappa_{i,j}^b+\kappa_{i,j}^t\big)u^\epsilon_{i,j}+\kappa_{i,j}^bu^\epsilon_{i,j-1}+f_{i,j}=\kappa_{i,j}^l\big(u^\epsilon_{i,j}-u^\epsilon_{i-1,j}\big)-\kappa_{i,j}^r\big(u^\epsilon_{i+1,j}-u^\epsilon_{i,j}\big).
\eeq
Then comparing the summation of the above equation for all $i\in\{1,\cdots, I-1\}$ with the boundary condition
\beq\label{eq:5dfmbound}
\sum_{i=1}^{I-1}\Big(\kappa_{i,j}^tu^\epsilon_{i,j+1}-(\kappa_{i,j}^t+\kappa_{i,j}^b)u^\epsilon_{i,j}+\kappa_{i,j}^bu^\epsilon_{i,j-1}\Big)
=-\sum_{i=1}^{I-1}f_{i,j},
\eeq
yields
$$
k^l_{1,j}(u^\epsilon_{1,j}-u^\epsilon_{0,j})=k^r_{I-1,j}(u^\epsilon_{I,j}-u^\epsilon_{I-1,j}).
$$
Using $u_{1,j}=u_{0,j}$ gives $u_{I,j}=u_{I-1,j}$. Therefore, the system produces the same results as the $5$-point FDM coupled with
 the  Neumann  boundary discretization such that $u^\epsilon_{1,j}=u^\epsilon_{0,j}$, $u^\epsilon_{I,j}=u^\epsilon_{I-1,j}$.

The advantage of using \eqref{eq:5dfmbound} is when $\epsilon\to 0$, in which case $\kappa_{i,j}^l\to +\infty$ and $\kappa_{i,j}^r\to +\infty$ for
$\forall i\in\{1,2,\cdots, I-1\}$ in \eqref{eq:5dfminternal}.
Together with $u_{1,j}=u_{0,j}$, one gets at the leading order \beq u^0_{0,j}=u^0_{1,j}=u^0_{2,j}=\cdots=u^0_{I,j}\equiv \bar{u}^{0}_{j}.\label{eq:limuj}\eeq
Then \eqref{eq:5dfmbound} gives the equation to determine the value of $\bar{u}^0_j$ such that
$$
\sum_{i=1}^{I-1}\kappa_{i,j}^t\bar{u}^0_{j+1}-\sum_{i=1}^{I-1}(\kappa_{i,j}^t+\kappa_{i,j}^b)\bar{u}^0_{j}+\sum_{i=1}^{I-1}\kappa_{i,j}^b\bar{u}^0_{j-1}
=-\sum_{i=1}^{I-1}f_{i,j},\qquad\mbox{for $j=1,2,\cdots, J$}.
$$
Therefore the solution is uniquely determined in the limit.
 \qed
 The convergence order can be improved by using a ghost point at the Neumann boundary and a better approximation of the integral in the third equation in \eqref{eq:ellipT3}.

{\begin{remark}
The condition number of the coefficient matrix for the discretized system can be considered as a continuous function of $\epsilon$. We denote it by $C(\epsilon)$. On the one hand, the limiting matrix is independent of $\epsilon$ and the solution can be uniquely determined	 from Theorem \ref{theo:unique}. This indicates that $C(0)$ is bounded. On the other hand, for fixed $\epsilon$, the discretized system is equivalent to the standard five point finite difference scheme coupled with Neumann boundary conditions on the left and right boundaries, thus $C(\epsilon)$ is bounded for fixed $\epsilon$. From the continuity of the coefficient matrix with respect to $\epsilon$, there exists a uniform bound of $C(\epsilon)$ for all $\epsilon\in[0,1]$.
\end{remark}
}

\subsection{$9$-Point FDM for the Nonaligned Case}
\begin{figure}[htb]
\centering
{
\includegraphics[width=7.5cm] {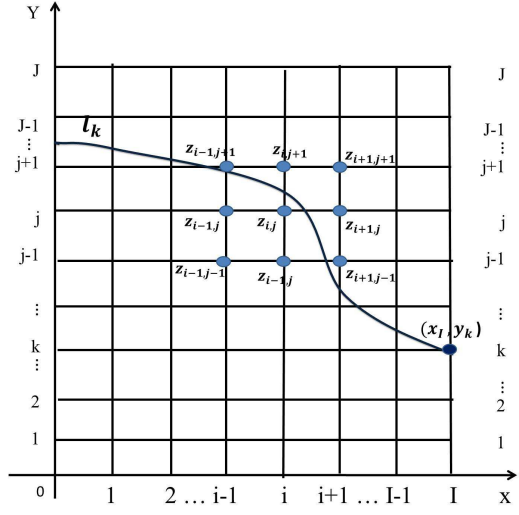}}
\caption{
{Notations for the discrete integral equation along $b$ field}.}\label{fig:9point}
\end{figure}

{When $\theta\neq 0$, we use the $9$-Point FDM to discretize the system \eqref{eq:SimellipC}. To simplify the notations as well as to make the discussions clearer, we consider the field line as in Figure
 \ref{fig:9point}.
We assume that $\Gamma_{in}$ is the left boundary and
$\Gamma_{out}$ is the right, while Dirichlet boundary conditions are imposed on the top and bottom.}

For all those grid points inside the computational domain, we use the classical $9$-Point FDM. The details are as follows: First of all, we introduce the notations such that $u_{i\pm 1/2,j}\approx u(x_{i\pm 1/2},y_j)$, $u_{i,j\pm 1/2}\approx u(x_{i},y_{j\pm 1/2})$ for any function $u$. The standard center difference approximations
read
 \beq\label{eq:9point}\begin{array}{ll}
\displaystyle \frac{\partial u^{\epsilon}}{\partial x}\Big|_{i\pm 1/2,j}\approx\frac{u_{i\pm1,j}^{\epsilon}-u_{i,j}^{\epsilon}}{\pm h_x},\quad \frac{\partial u^{\epsilon}}{\partial y}\Big|_{i,j\pm 1/2}\approx\frac{u_{i,j\pm 1}^{\epsilon}-u_{i,j}^{\epsilon}}{\pm h_y},\\
\displaystyle \frac{\partial u^{\epsilon}}{\partial y}\Big|_{i\pm 1/2,j}\approx\frac{u_{i\pm 1/2,j+1}^{\epsilon}-u_{i\pm 1/2,j-1}^{\epsilon}}{2h_y}\approx\frac{u_{i\pm 1,j+1}^{\epsilon}+u_{i,j+1}^{\epsilon}-u_{i,j-1}^{\epsilon}
-u_{i\pm 1,j-1}^{\epsilon}}{4h_y},\\
\displaystyle \frac{\partial u^{\epsilon}}{\partial x}\Big|_{i,j\pm 1/2}\approx\frac{u_{i+1,j\pm 1/2}^{\epsilon}-u_{i-1,j\pm 1/2}^{\epsilon}}{2h_y}\approx\frac{u_{i+1,j\pm 1}^{\epsilon}+u_{i+1,j}^{\epsilon}-u_{i-1,j}^{\epsilon}-u_{i-1,j\pm 1}^{\epsilon}}{4h_x}.
\end{array}
\eeq
 Let $Q=A\nabla u^{\epsilon}$, we can approximate
the diffusion operator at $(x_i,y_j)$ by
 \beq\label{eq:9point2}
\nabla \cdot (A\nabla u^{\epsilon})|_{i,j}=\f{\p Q}{\p x}\Big|_{i,j}+
\f{\p Q}{\p y}\Big|_{i,j}\approx\frac{Q_{i+1/2,j}-Q_{i-1/2,j}}{h_x}+\frac{Q_{i,j+1/2}-Q_{i,j-1/2}}{h_y}.
\eeq
with
 $$
Q_{i\pm 1/2,j} \approx A_{i\pm 1/2,j}\cdot\left(\frac{\partial u^{\epsilon}}{\partial x}\Big|_{i\pm 1/2,j},\frac{\partial u^{\epsilon}}{\partial y}\Big|_{i\pm 1/2,j}\right)^T,\qquad
Q_{i,j\pm 1/2} \approx A_{i,j\pm 1/2}\cdot\left(\frac{\partial u^{\epsilon}}{\partial x}\Big|_{i,j\pm 1/2},\frac{\partial u^{\epsilon}}{\partial y}\Big|_{i,j\pm 1/2}\right)^T.
 $$

The Neumann boundary conditions at $\Gamma_{in}$ are discretized locally by using
\beq\label{eq:nbound}\frac{\partial u^{\epsilon}}{\partial y}\Big|_{0,j}\approx\frac{u_{0,j+1}^{\epsilon}-u_{0,j-1}^{\epsilon}}{2 h_y},\quad \frac{\partial u^{\epsilon}}{\partial x}\Big|_{0,j}\approx\frac{-\f{3}{2} u_{0,j}^{\epsilon}+2u_{1,j}^{\epsilon}-\f{1}{2} u_{2,j}^{\epsilon}}{h_x}.\eeq
We let $u_{i,j}=0, (i=0,I,j=0,1,\cdots J)$ to satisfy the Dirichlet boundary  conditions on $\Gamma_D$.

The major difference and difficulty is about the discretization of the integration in \eqref{eq:SimellipC}. The details are as follows:
\begin{itemize}
\item  The first step is to determine those field lines and intersection points for integration. To simplify the notation, we assume that all field lines start at the left side and end at the right side of $\Omega$. In order to have the number of equations being the same as the number of unknowns,
it is important to choose those field lines that pass through those grid points on the right. Assume $(x_I, y_k)\in \Gamma_{out}$ , we define $l_k$ as the field line that passes
through $(x_I, y_k)$, $(k=1\cdots J-1)$.

Based on the vector field $\mathbf{b}=(\cos\theta,\sin\theta)$, we can determine the field line by the following nonlinear ODE,
 \beq\label{eq:ODE}
\left\{\begin{array}{ll}
\displaystyle \dot x(t)=\cos\theta(x(t),y(t))\quad x(0)=x_I,
\\
\displaystyle \dot y(t)=\sin\theta(x(t),y(t))\quad y(0)=y_k.
\end{array}
\right.
\eeq
 To approximate the field line,
we use second order Runge-Kutta Method with a fine mesh to solve the above system. In order to get the integration along the field line $l_k$, we have to determine the quadrature points that are given by the intersection points of
 $l_k$ with those cell edges that are parallel to the $y$ axis. We denote the intersection points by $(x_{i},\bar{y}^k_{i})$,
for $i\in\{0,\cdots,I\}$, and we assume $\bar{y}^k_{i}\in[y_{k_i},y_{k_i}+h_y)$ with $k_i\in\{0,\cdots,J-1\}$.

     \item Approximate $E\mid_{(x_{i},\bar{y}^k_{i})}=\exp\Big(\int_0^ {(x_{i},\bar{y}^k_{i})} \nabla \cdot \mathbf b^T  \,\,ds\Big)$ by trapezoidal's rule.
     $E\mid_{(x_{0},\bar{y}^{k}_{0})}=1$ and when $i>0$,
     \beq\label{eq:Eij}
{\int_0^ {(x_i,\bar{y}^k_{i})} \nabla \cdot \mathbf b^T  \,\,ds}\approx\displaystyle \sum_{m=1}^{{i-1}}\nabla\cdot\mathbf{b}^T\mid_{(x_{m},\bar{y}^k_{m})}\frac{h_x}{\cos\theta(x_{m},\bar{y}^k_{m})}+\f{1}{2}\sum_{m=0,i} \nabla\cdot\mathbf{b}^T\mid_{(x_{m},\bar{y}^k_{m})}\frac{h_x}{\cos\theta(x_{m}.\bar{y}^k_{m})},
\eeq
     \item The diffusion operator $- \nabla
   \cdot \big( \alpha \mathbf{b}_{\perp} (\mathbf{b}_{\perp}\cdot \nabla u^{\epsilon})\big)\big|_{i,j}$ is approximated by centered finite difference method similar as in \eqref{eq:9point2}. The value of $\Theta=f+\nabla
   \cdot \big( \alpha \mathbf{b}_{\perp} (\mathbf{b}_{\perp}\cdot \nabla u^{\epsilon})\big)$ at the integration point
   $(x_{i},\bar y^k_{i})$ can be given by linear interpolation such that:
   \beq
   \begin{array}{ll}
  \Theta\mid_{(x_{i},\bar y^k_{i})}\approx \Theta\mid_{(x_{i},y_{k_i})}\frac{h_{d_{i}}}{h_y}+\Theta\mid_{(x_{i},y_{k_i}+h_y)}\frac{h_y-h_{d_{i}}}{h_y},
    \end{array}
  \eeq
  where $h_{d_{i}}=y^{k}_{i}+h_y-\bar y^{k}_{i}$.
  \item
 For each $l_k$, $k=1,2,\cdots,J-1$, trapezoidal's rule is employed to discretize $\int_0^{L_l} E\big( \nabla \cdot ( \alpha
   \mathbf b_{\perp} \mathbf b_{\perp}^T \nabla u^{\epsilon}) +f\big){ds}$.
    The third equation of \eqref{eq:SimellipC} can be approximated by
\beq\label{eq:discL}\begin{aligned} &\sum_{i=1}^{I-1}E\mid_{(x_i,\bar{y}^k_{i})}\Theta\mid_{(x_{i},\bar{y}^k_{i})}\times\frac{h_x}{\cos\theta(x_{i},\bar{y}^k_{i})}
  +\sum_{i=0,I}E\mid_{(x_{i},\bar{y}^k_{i})} \Theta\mid_{(x_{i},\bar{y}^k_{i})}
     \times\frac{h_x/2}{\cos\theta(x_{i},\bar{y}^k_{i})}
     \\&
    =\alpha E\f{\mathbf{n}\cdot\mathbf{b}^T_\bot}
    {\mathbf{n}\cdot\mathbf{b}^T}\big(-\sin\theta\p_xu^\epsilon+\cos\theta\p_yu^\epsilon\big)\Big|_{(x_{0},\bar{y}^k_{0})}^{(x_{I},\bar{y}^k_{I})}
         \end{aligned}\eeq
         where $\p_xu^\epsilon\mid_{(x_{i},y^k_{i})}$ and $\p_yu^\epsilon\mid_{(x_{i},y^k_{i})}$ $(i=0, I)$ are discretized similar as \eqref{eq:nbound}.
\end{itemize}


\begin{remark} Here the easiest and most straight forward FDM is used, but other discretizations can be applied as well. If there exits grid point $(x_i,y_0)$ at the bottom side of $\Omega$ that belongs to $\Gamma_{out}$, we only need to replace the Neumann boundary condition at $(x_i,y_0)$ by the integration over the field line that passes through $(x_i,y_0)$. However to get a good approximation of the field line integration, the intersection points with some edges parallel to the $x$ axis might be needed.
    We assume that the field lines start on the left side of $\Omega$ and end on the right side of $\Omega$ is to simplify the notations.
\end{remark}


\section{Numerical Results}\label{num}
We present several computational tests to demonstrate the performance of our proposed scheme. The first tree examples are for the aligned case and the last two are for the non-aligned non-uniform $\mathbf{b}$ field.


\textbf{Example 1:} Space uniform and aligned case. Let the field $\mathbf{b}$ be aligned with $x$ axis as in \eqref{eq:ellipT3}.
 The computational domain is $\Omega=[0,1]\times[0,1]$ and $\alpha\equiv 1$
is space uniform.
The source term $f$ is given by
$$
f=(4+\epsilon)\pi^2\cos(2\pi x)\sin(\pi y)+\pi^2\sin(\pi y).
$$
Then the exact solution $u_{exact}$ of \eqref{eq:ellipT1} is
$$
u_{exact}=\sin(\pi y)+\epsilon \cos(2\pi x)\sin(\pi y).
$$

\begin{figure}[htb]
\centering
{
\subfigure [] {\includegraphics[width=7.5cm,clip]{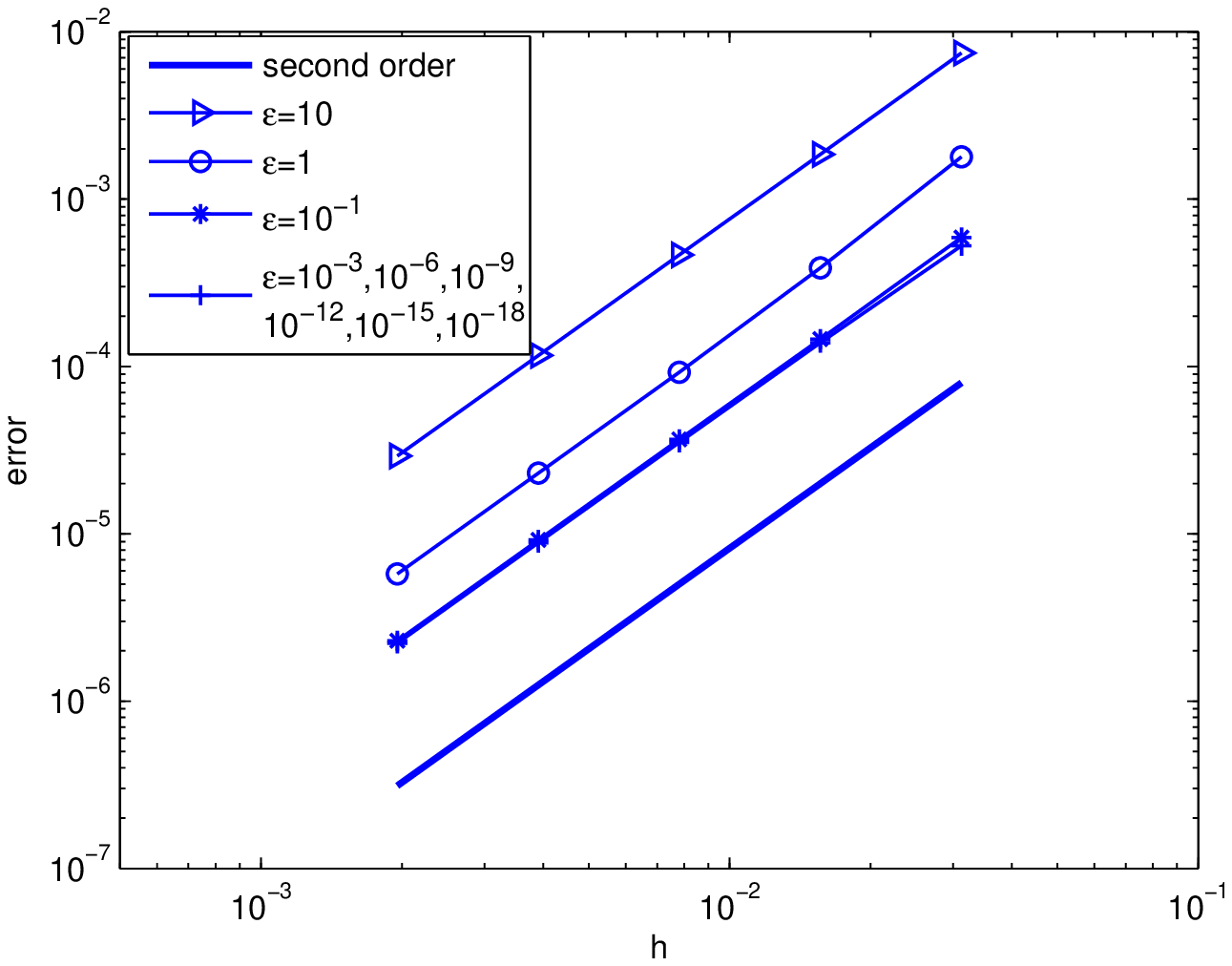} }
\subfigure [] {\includegraphics[width=7.5cm,clip]{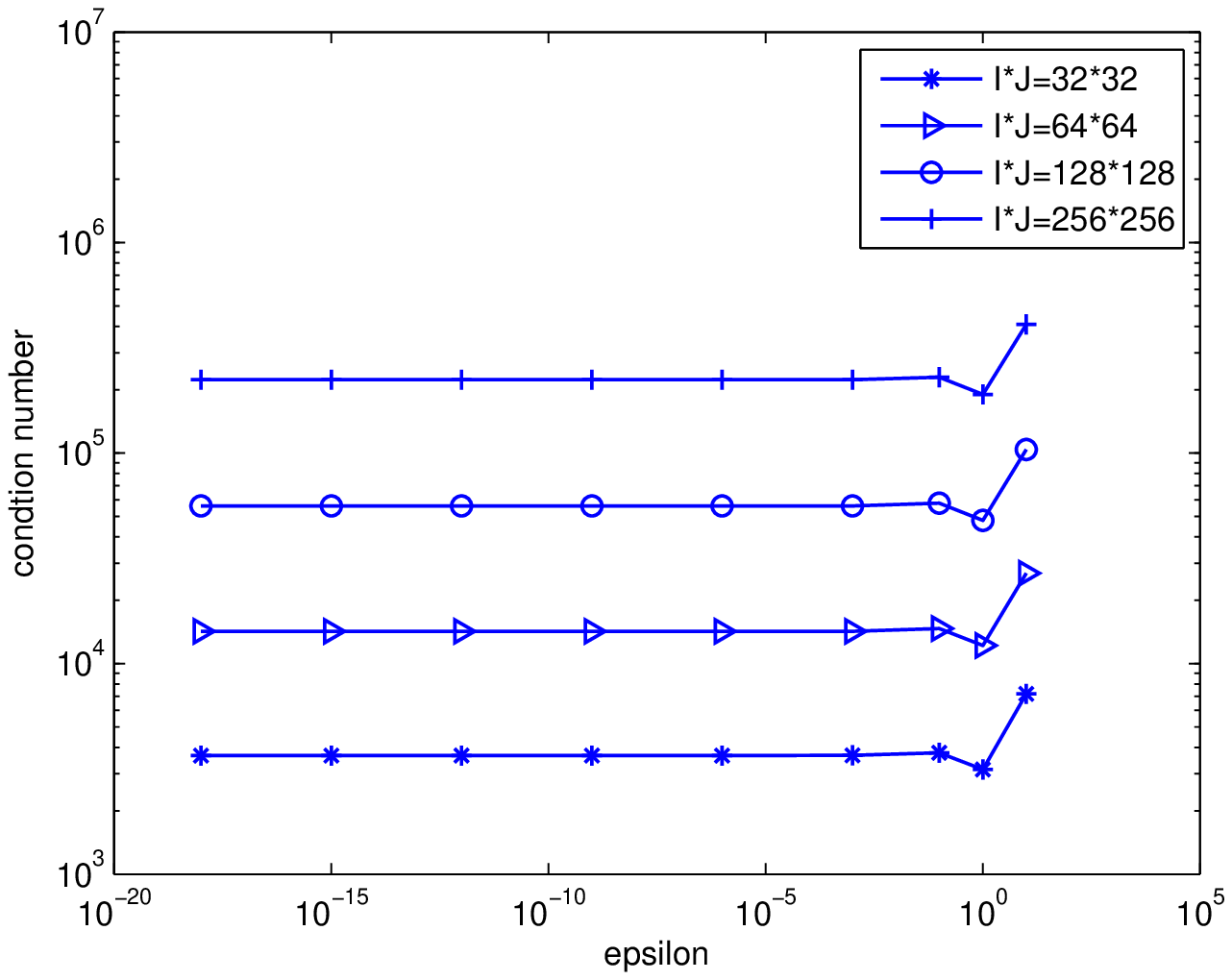} }

}
\caption{Example 1.  The numerical convergence order and condition number of the coefficient matrices for the $5$-point FDM. The results of different grids and different $\epsilon$ are given. (a) Convergence orders; (b) Condition number}.\label{fig:exone}
\end{figure}

For different values of $\epsilon$ that ranges from $10^{-18}$ to $10$, we present the discrete $L^2$ norms of the errors and their corresponding convergence orders in Table \ref{tab:1} and  Figure \ref{fig:exone} $(a)$.
 Uniform second order convergence can be observed. Furthermore, as plotted in Figure \ref{fig:exone} $(b)$, the condition number of the new scheme is bounded by
 a constant independent of $\epsilon$, whose magnitude coincides with the AP Micro-Macro decomposition scheme in \cite{Degond121} and the Duality-Based scheme in \cite{Degond122}.

 \begin{table}[tbp]
\centering
 \begin{tabular}{|cccccc|}
\hline
$ \epsilon\backslash I\times J$  &$32\times32$  & $64\times64$ & $128\times128$ & $256\times256$& $512\times512$   \\ \hline
10 & 7.4692E-03  & 1.8565E-03   &4.6476E-4  & 1.1645E-04  & 2.9160E-05 \\ \hline
1  & 1.7929E-03   & 3.8911E-04   &9.2325E-05  & 2.3095E-05  & 5.7580E-06\\ \hline
$10^{-1}$ & 5.8925E-04  & 1.4572E-04    &3.6610E-05  &9.1794E-06 & 2.2985E-06 \\ \hline
$10^{-3}$  & 5.2710E-04  & 1.3883E-04   &3.5228E-05  & 8.8486E-06  & 2.2158E-06\\ \hline
$10^{-6}$  & 5.2710E-04  & 1.3883E-04   &3.5228E-05  & 8.8486E-06  & 2.2158E-06\\ \hline
$10^{-9}$  & 5.2710E-04 & 1.3883E-04   &3.5228E-05  & 8.8486E-06  & 2.2158E-06\\ \hline
$10^{-12}$  & 5.2710E-04  & 1.3883E-04   &3.5228E-05  & 8.8486E-06  & 2.2158E-06\\ \hline
$10^{-15}$  & 5.2710E-04  & 1.3883E-04   &3.5228E-05  & 8.8486E-06  & 2.2158E-06\\ \hline
$10^{-18}$   & 5.2710E-04  & 1.3883E-04   &3.5228E-05  & 8.8486E-06  & 2.2158E-06\\ \hline
\end{tabular}
 \caption{Example 1. The discrete $L^2$ norm of the errors for different grids and $\epsilon$ values.}
\label{tab:1}
\end{table}

%
%
\bigskip

\begin{figure}[htb]
\centering
{
\subfigure [] {\includegraphics[width=7.5cm,clip]{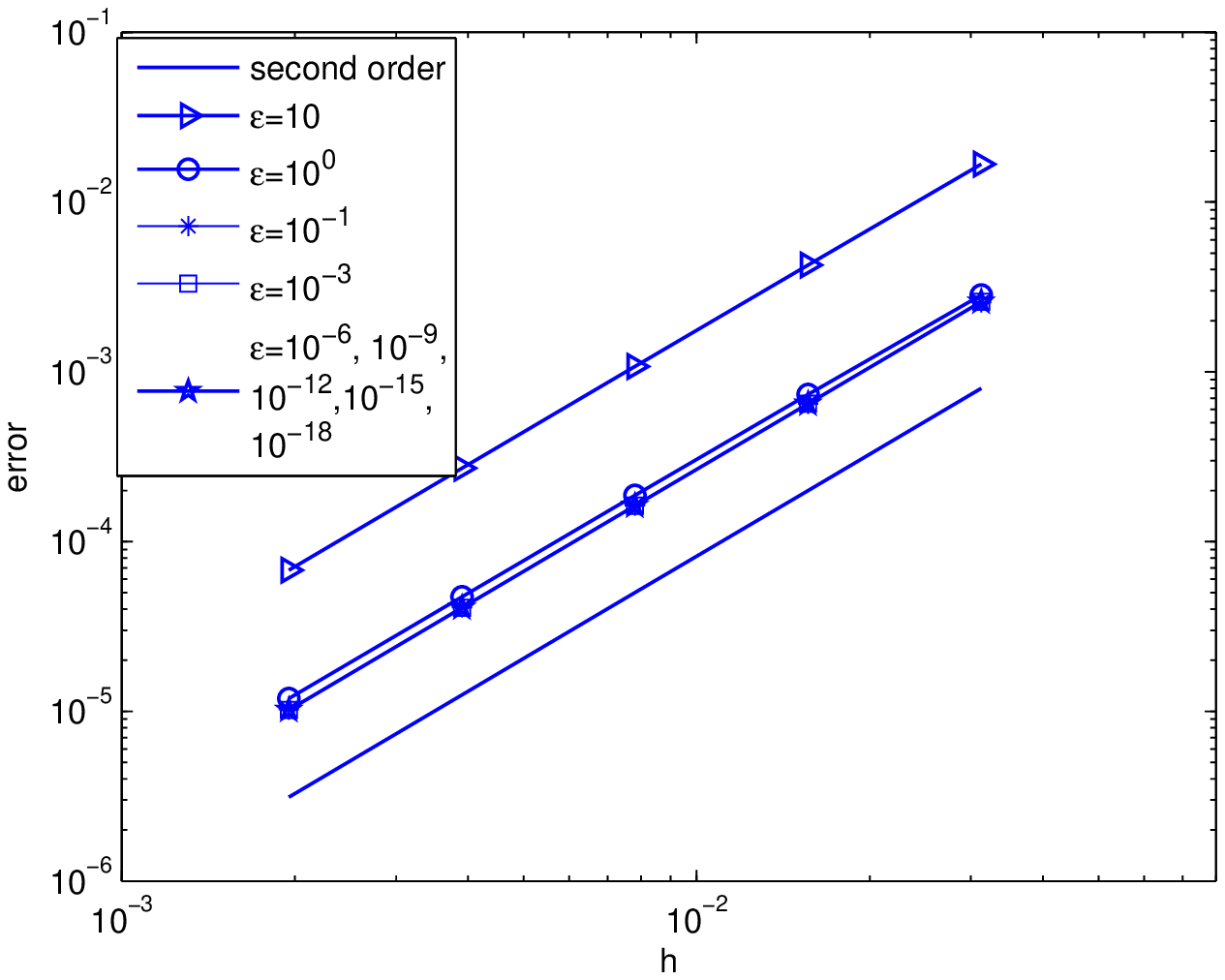} }
\subfigure [] {\includegraphics[width=7.5cm,clip]{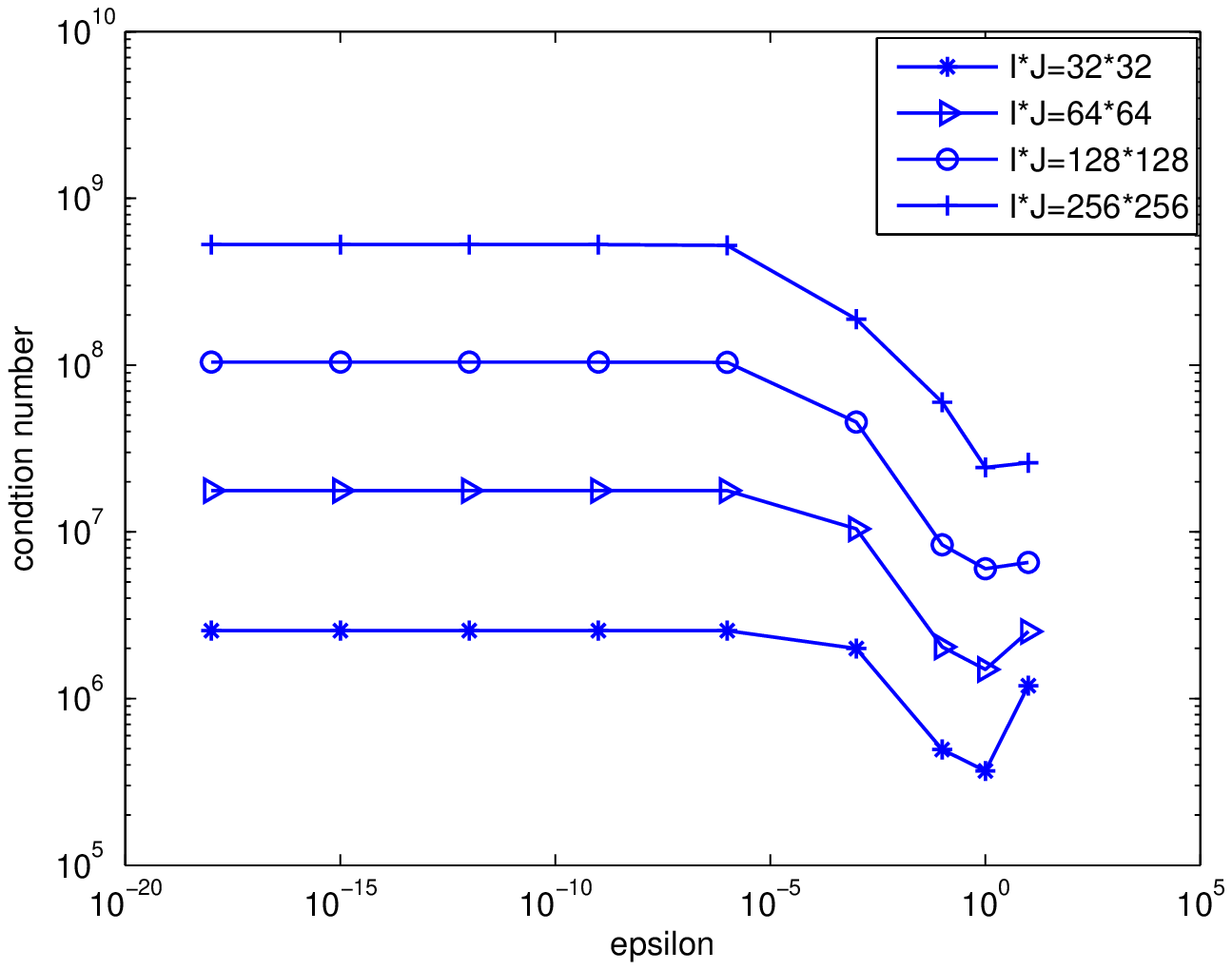} }

}
\caption{Example 2. The numerical convergence order and condition number of the coefficient matrices for the $5$-point FDM.
The results with different grids for different $\epsilon$ are given. (a) Convergence orders; (b) Condition number.}\label{fig:extwo}
\end{figure}

\textbf{Example 2:} Space non-uniform and aligned case. We choose  \beq\label{eq:test2A}
 A(x,y)=\left(\begin{array}{cc}c_1+xy^2&0\\0&\frac{1}{\epsilon}(c_2+xy)\end{array}\right),
\eeq
with $c_1$, $c_2$ two positive constants. Let the exact solution be
 \beq\label{eq:testtwo}
u_{exact}=\sin(\frac{2\pi}{L_x}x)+\epsilon\cos(\frac{2\pi}{L_y}y)\sin(\frac{2\pi}{L_x}x).
\eeq
We choose $L_x=L_y=10$ and $c_1=c_2=L_y$ and the computational domain is $[0,L_x]\times[0,L_y]$. The right-hand source term $f$ is determined by inserting \eqref{eq:testtwo} into  \eqref{eq:ellipT1}.  This example is the same as the one in \cite{Degond101}, where lagrange multiplier is used to maintain the well-posedness in the limit. The computational cost our proposed approach is less expensive. As can be seen from Table \ref{tab:2} and Figure \ref{fig:extwo}, the new scheme is of second order accuracy and the condition number of the discretized system is bounded by a constant independent of $\epsilon$.

 \begin{table}[tbp]
\centering
 \begin{tabular}{|cccccc|}
\hline
$ \epsilon\backslash I\times J$  &$32\times32$  & $64\times64$ & $128\times128$ & $256\times256$& $512\times512$   \\ \hline
10 &1.6751E-02	&4.2638E-03	&1.0767E-03	&2.7064E-04	&6.7851E-05
 \\ \hline
1          &2.8313E-03	&7.3399E-04	&1.8713E-04	&4.7267E-05	&1.1879E-05
\\ \hline
$10^{-1}$ & 2.5698E-03	&6.4871E-04	&1.6317E-04	&4.0935E-05	&1.0253E-05
\\ \hline
$10^{-3}$  &2.5866E-03	&6.4852E-04	&1.6253E-04	&4.0697E-05 & 1.0167E-05\\ \hline
$10^{-6}$  	&2.5881E-03	&6.4905E-04	&1.6269E-04	&4.0740E-05&1.0186E-05\\ \hline
$10^{-9}$  & 2.5881E-03	&6.4906E-04	&1.6269E-04	&4.0740E-05  &1.0186E-05\\ \hline
$10^{-12}$  & 2.5881E-03	&6.4906E-04	&1.6269E-04	&4.0740E-05 &1.0186E-05\\ \hline
$10^{-15}$ & 2.5881E-03	&6.4906E-04	&1.6269E-04	&4.0740E-05 &1.0186E-05\\ \hline
$10^{-18}$  & 2.5881E-03	&6.4906E-04	&1.6269E-04	&4.0740E-05 &1.0186E-05\\ \hline
\end{tabular}
 \caption{Example 2. the discrete $L^2$ norm of the errors for different grids and $\epsilon$ values.}
\label{tab:2}
\end{table}

\bigskip

\textbf{Example 3:} Variable $\epsilon$ and {aligned} $\mathbf{b}$ field case. We test a case when $\epsilon$ varies from $1$ to some fixed parameter $\epsilon_{\min}$ in this example.
 Let
\beq\label{eq:vepsilon}
\epsilon(x,y)=\f{1}{2}[1+\tanh (a(x_0-x))+\epsilon_{\min}(1-\tanh a(x_0-x))],\qquad \alpha\equiv 1.
\eeq
Here {$a$} and $x_0$ respectively control the width and position of the transition region. In our simulations,
we set $x_0=0.25$ and $a=50$. The anisotropy varies smoothly in the whole computational domain, and changes its value in a relatively narrow transition region.

\begin{figure}[htb]
\centering
{
\subfigure [] {\includegraphics[width=7.5cm,clip]{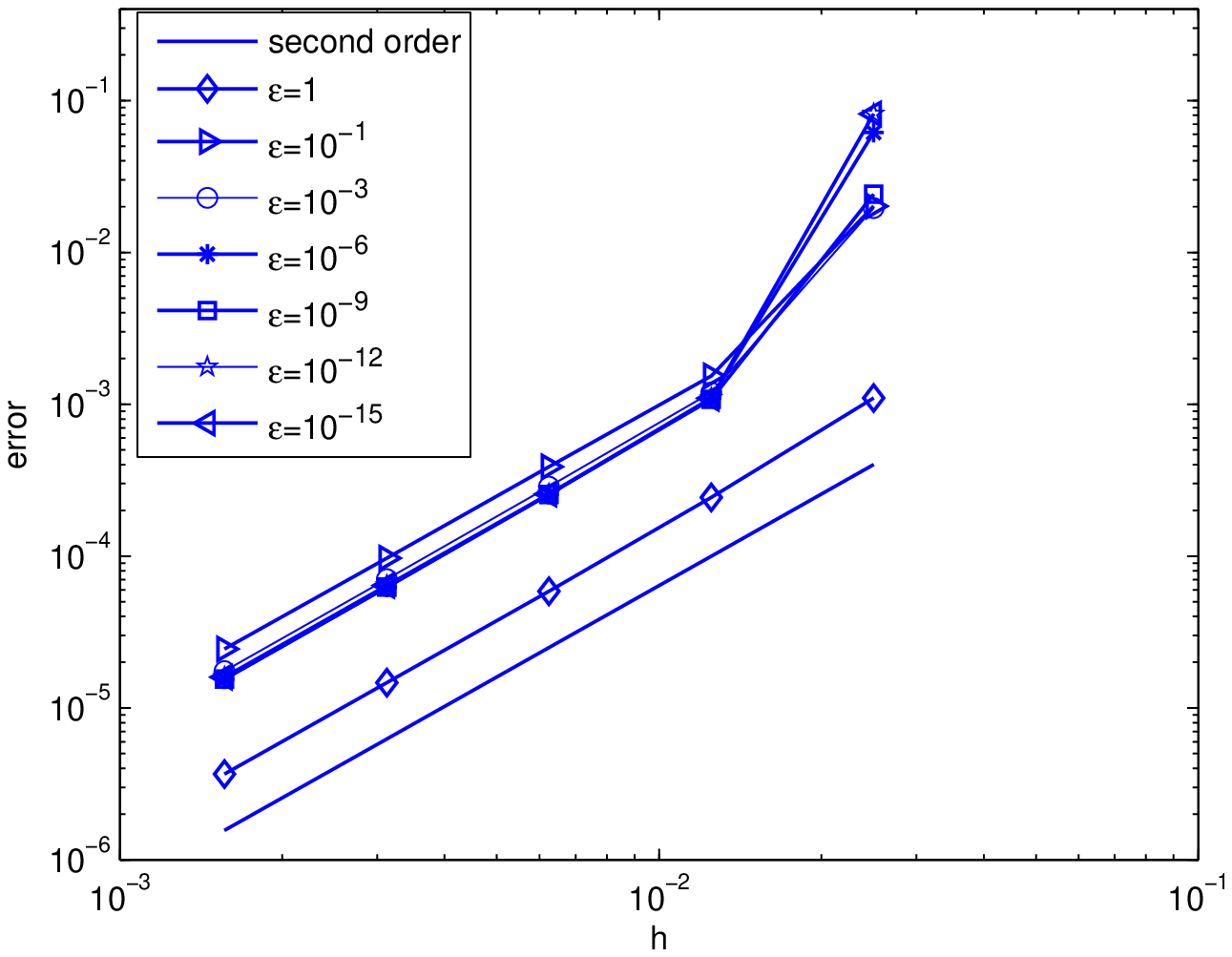} }
\subfigure [] {\includegraphics[width=7.5cm,clip]{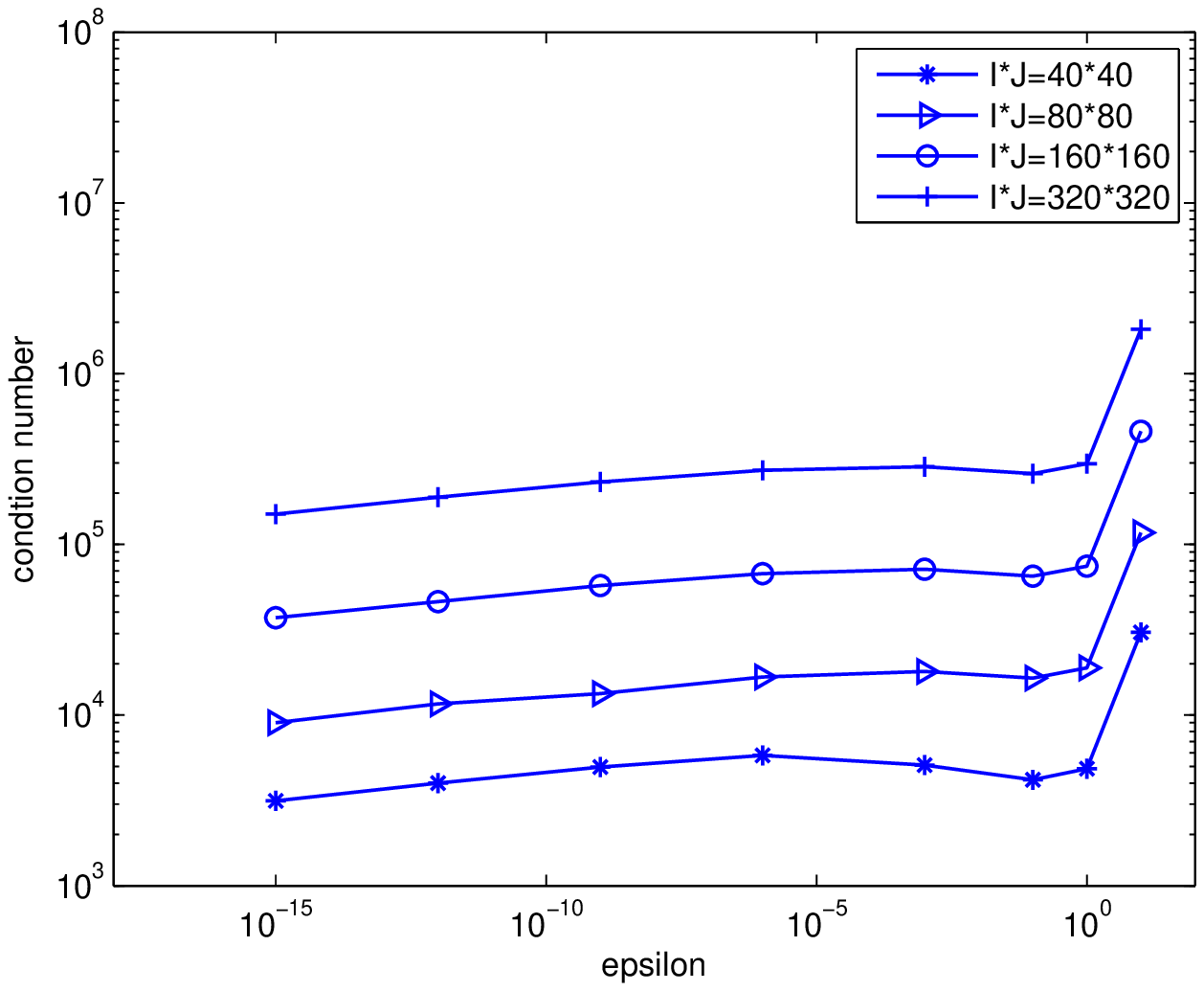} }

}
\caption{Example 3. Convergence orders and Condition number estimate for the discretization matrices of the $5$-point FDM. Different grids and different $\epsilon$ values are used. (a) Convergence orders; (b) Condition number.}\label{fig:exthree}
\end{figure}

In \cite {Degond122}, the authors have pointed out that it is necessary to use an AP scheme when $\epsilon_{min}$ become small. Direct simulation of the original system \eqref{eq:ellipT1} yields an error increases with $1/\epsilon_{min}$.
Our scheme can achieve second order convergence for all $\epsilon_{\min}$ ranging from $10^{-15}$ to $10$. %
The results are displayed in Table \ref{tab:3} and Figure \ref{fig:exthree}. The new scheme is capable of producing accurate results when the strength of anisotropy varies
 a lot in the computational domain. When $\epsilon_{min}$ becomes small, $u^\epsilon$ is almost constant along the $x$ axis
 at the part $\epsilon(x,y)$ is small and the largest error appear at the
 transition part. {As in Table \ref{tab:3} and Figure \ref{fig:exthree},
 when we refine the mesh from $I\times J=40\times 40$
 to $I\times J=80\times 80$, the errors decrease faster than second order. This is because
 $a=50$ in our simulations which yield a narrow transition region. The mesh size of $I\times J=40\times 40$ is too
 large to accurately represent the narrow transition region, which, therefore, introduces a relatively large error.
 The effect of the narrow transition region becomes significant when $\epsilon\leq 10^{-1}$, which explains our numerical observations.}
 \begin{table}[tbp]
\centering
 \begin{tabular}{|cccccc|}
\hline
$ \epsilon\backslash I\times J$  &$40\times40$  & $80\times80$ & $160\times160$ & $320\times320$& $640\times640$   \\ \hline
10 & 4.8524E-02	&7.8469E-03	&1.9855E-03	&4.9790E-04 & 1.2454E-04
\\ \hline
1  &1.0987E-03	&2.4321E-04	&5.8752E-05	&1.4648E-05 &3.6737E-06
\\ \hline
$10^{-1}$ &2.0177E-02	&1.5356E-03	&3.8728E-04	&9.6871E-05 &2.4426E-05
 \\ \hline
$10^{-3}$  & 1.9554E-02	&1.1952E-03	&2.8670E-04	&7.0582E-05 &1.7465E-05
\\ \hline
$10^{-6}$ & 6.1764E-02	&1.0941E-03	&2.5508E-04	&6.2449E-05 &1.5503E-05

\\ \hline
$10^{-9}$  & 2.4117E-02	&1.0772E-03	&2.5495E-04	&6.2416E-05 &1.5462E-05

\\ \hline
$10^{-12}$  & 8.1520E-02	&1.0938E-03	&2.5608E-04	&6.4078E-05 &1.5963E-05

\\ \hline
$10^{-15}$ &2.0985E-02	&1.0787E-03	&2.5495E-04	&6.2415E-05 &1.5460E-05


\\ \hline
\end{tabular}
 \caption{Example 3. the discrete $L^2$ norm of the errors for different grids and $\epsilon$ values.}
\label{tab:3}
\end{table}

%

\bigskip

\textbf{Example 4:} {Uniform $\epsilon$ with non-aligned $\mathbf{b}$ field case.} We follow the construction of a test case in  \cite{Degond121,Degond122}
 that the exact solution can be analytically given, {but we vary the computational domain to include case when $\mathbf n \cdot
  \mathbf{b}^T_{\perp} \neq 0$ at $\Gamma_{out}$.}

\begin{figure}[htb]
\centering
{
\subfigure []  {\includegraphics[width=7.5cm,clip]{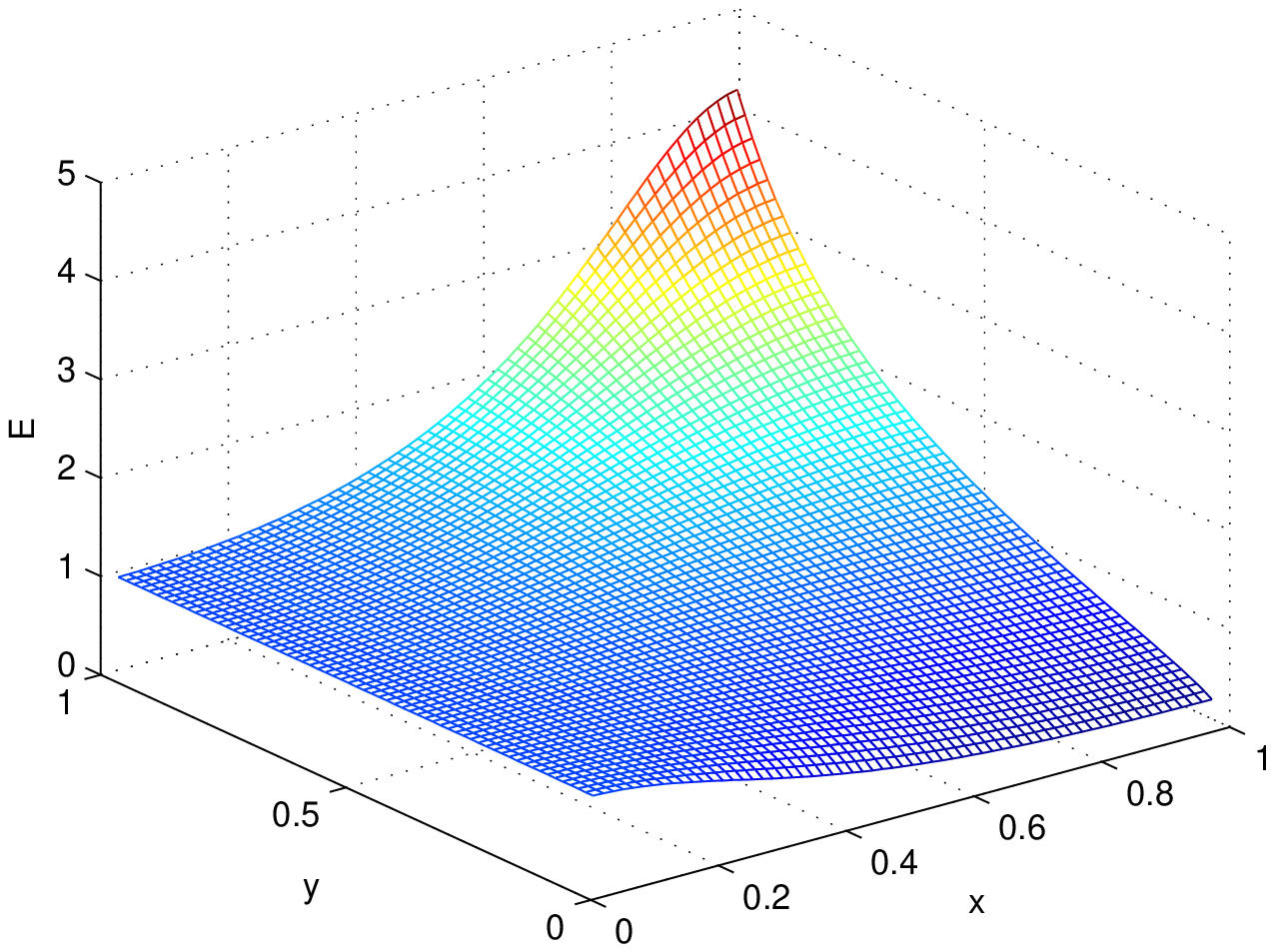} }
\subfigure [] {\includegraphics[width=7.5cm,clip]{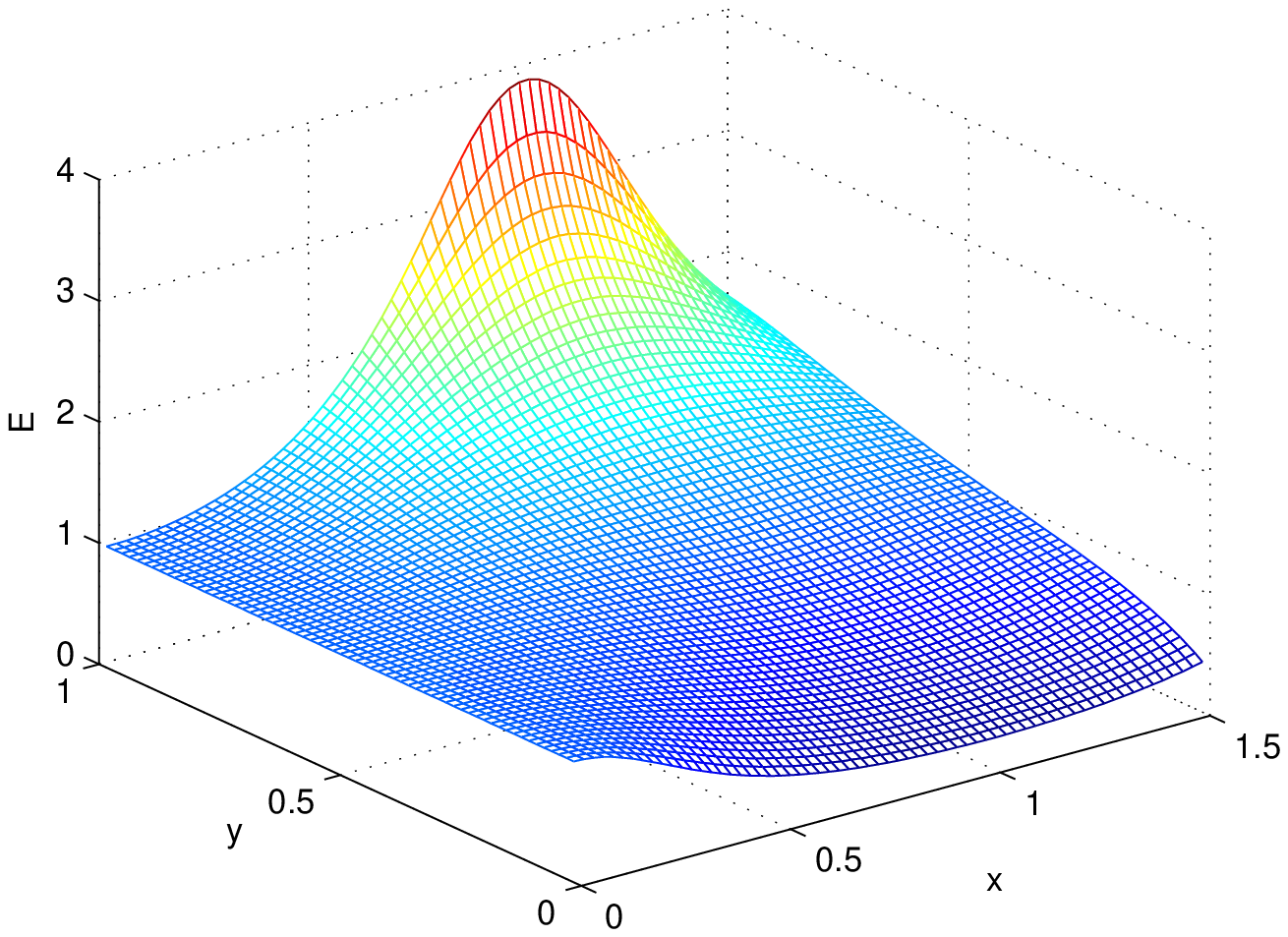} }

}
\caption{Example 4.
The numerical results for $E$ for different computational domain. (a) $[0,1]\times [0,1]$; (b)$[0,3/2]\times [0,1]$.}
\label{fig:exfourE}
\end{figure}

\begin{figure}[htb]
\centering
{
\subfigure [] {\includegraphics[width=7.5cm,clip]{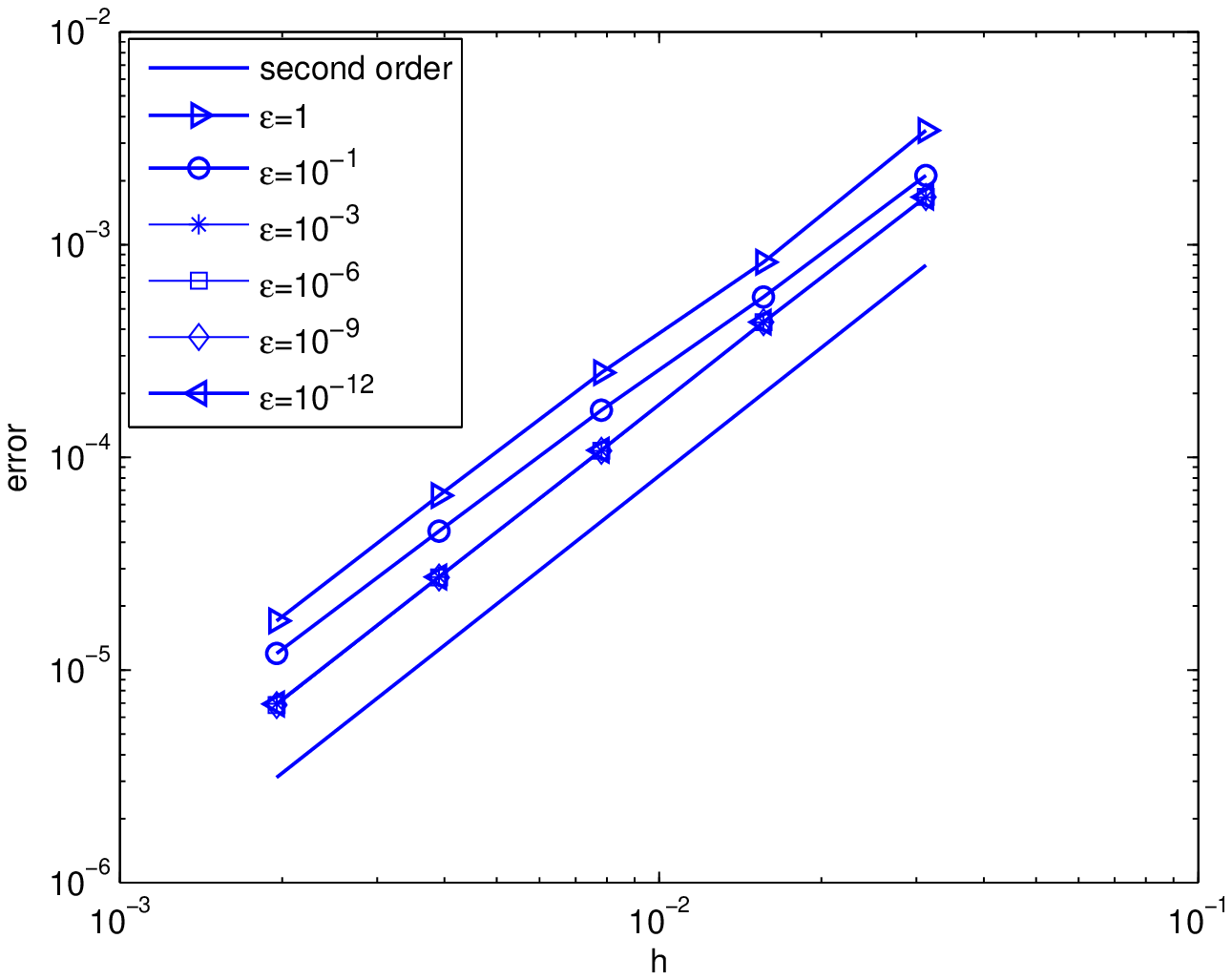} }
\subfigure [] {\includegraphics[width=7.5cm,clip]{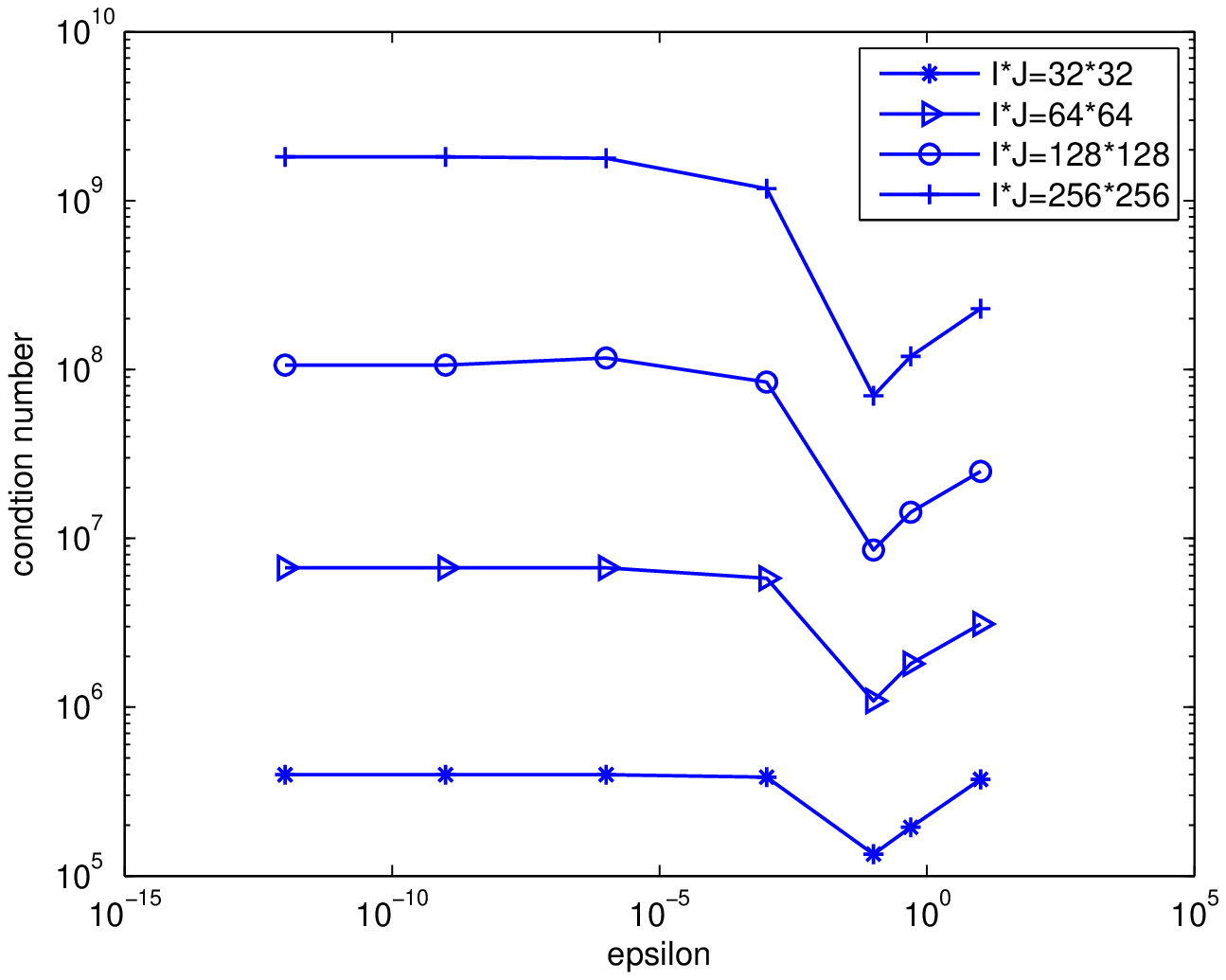} }

}
\caption{Example 4.
The numerical convergence order and condition number of the coefficient matrices for the $9$-point FDM.
The results with different grids for different $\epsilon$ are given. (a) Convergence orders; (b)Condition number.}\label{fig:exfour}
\end{figure}

\begin{figure}[htb]
\centering
{
\subfigure [] {\includegraphics[width=7.5cm,clip]{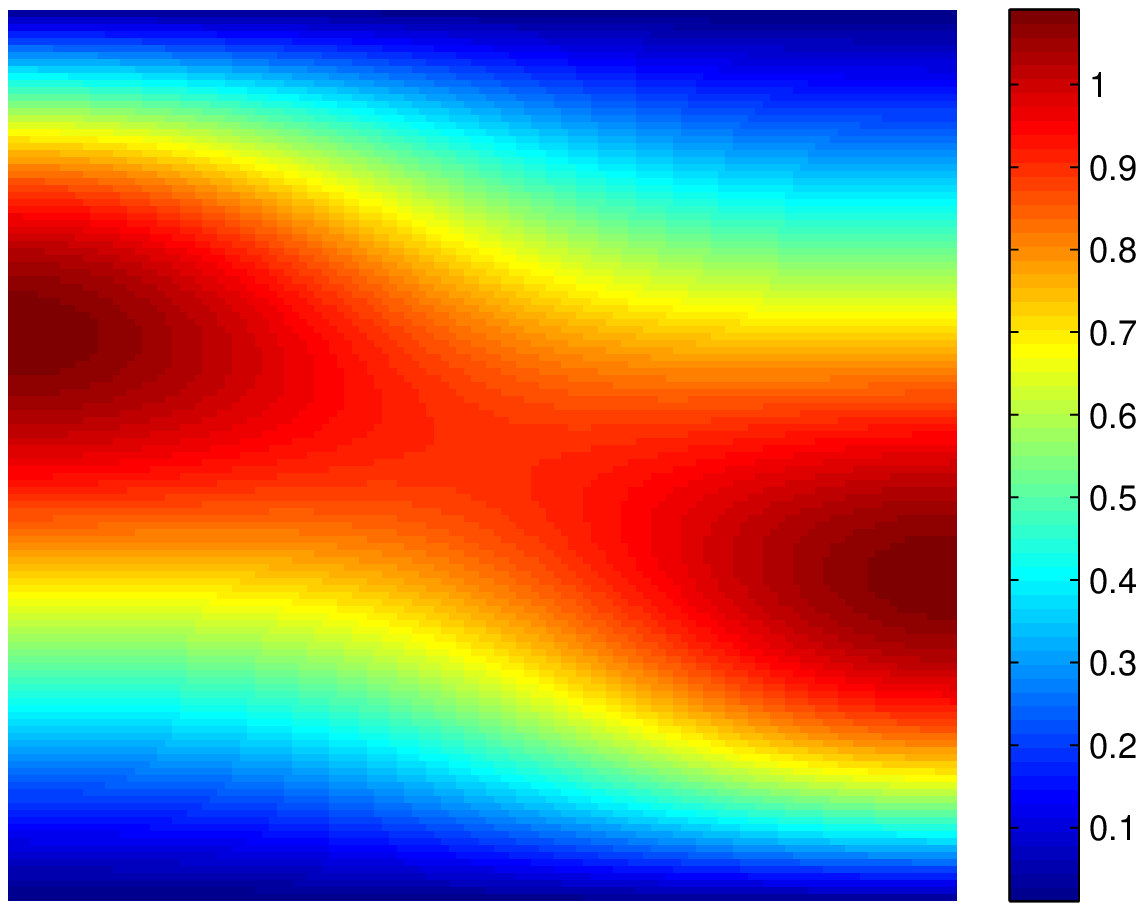} }
\subfigure[] {\includegraphics[width=7.5cm,clip]{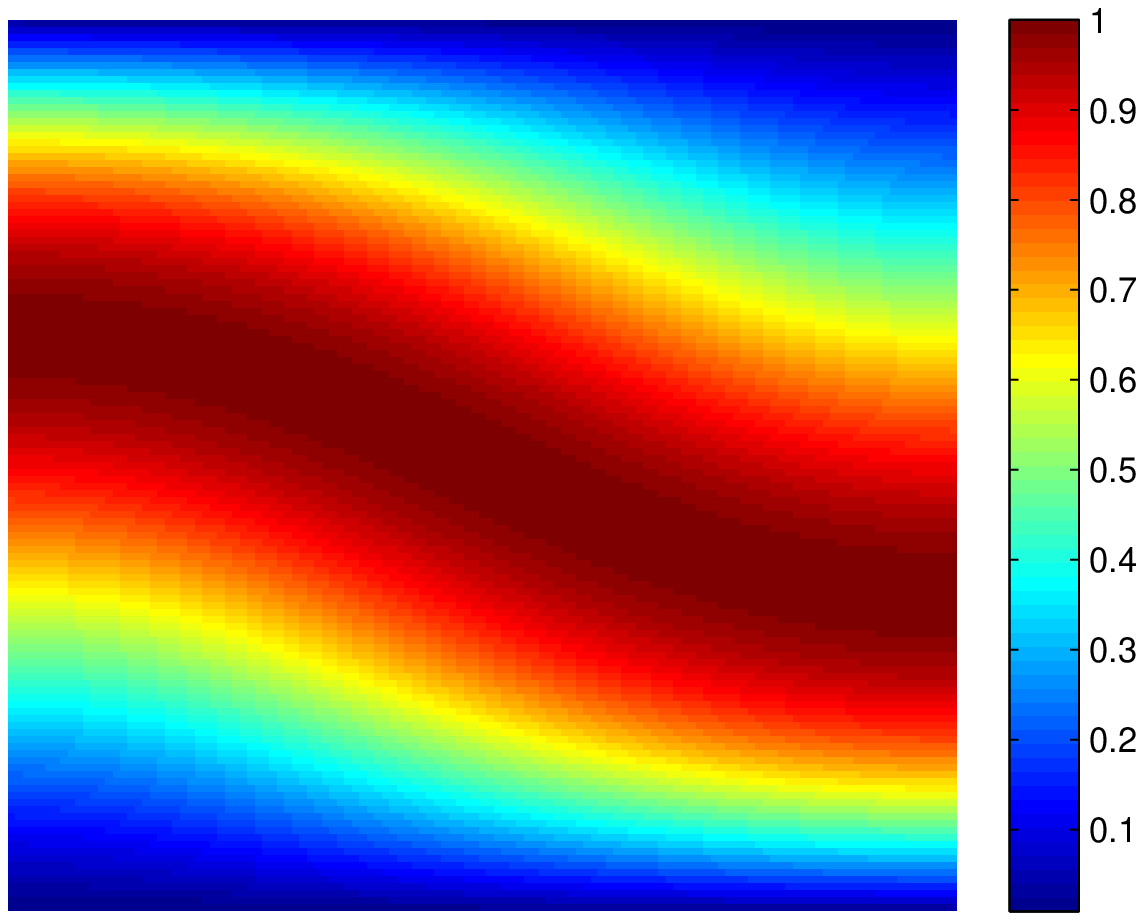}  }

}
\caption{Example 4.
 (a): $u^\epsilon(x,y)$ with $\epsilon=0.1$, $I\times J=128\times128$;
 (b): $u^\epsilon(x,y)$ with $\epsilon=10^{-12}$, $I\times J=128\times 128$.}
\label{fig:exfour2}
\end{figure}

 {\textbf{Test one:}}  Let the computational domain be $[0,1]\times[0,1]$ and the exact solution be
 \beq\label{eq:test3}
u^\epsilon_{exact}=\sin(\pi y+2(y^2-y)\cos(\pi x))+\epsilon\cos(2\pi x)\sin(\pi y).
\eeq
As $\epsilon \rightarrow 0$, $u^\epsilon_{excat}$ converges to $u^0=\sin(\pi y+2(y^2-y)\cos(\pi x))$.
Note that the field $\mathbf{b}$, satisfies $\mathbf{b}\cdot\nabla u^0=0$ which indicates that the limit solution $u^0$ is a constant along the $\mathbf{b}$ field line. This is
how the field line is constructed,
whose direction is given by
 \beq\label{eq:test3A}
\mathbf{b}=\frac{B}{\mid B\mid},\quad B=\left(\begin{array}{c}2 (2y-1)\cos(\pi x)+\pi\\2\pi (y^2-y)\sin(\pi x)\end{array}\right).
\eeq
The source term is calculated by plugging \eqref{eq:test3} into \eqref{eq:SimellipA}. Once again, as pointed out in \cite{Degond121,Degond122}, for fixed grids, direct simulation of the original system \eqref{eq:ellipT} yields an error increases with $1/\epsilon$ when $\epsilon$ is less than the order of $10^{-2}$. The $L^2$ and $H^1$ errors become $O(1)$ when $\epsilon$ is at $O(10^{-6})$ to $O(10^{-12})$, depending on the mesh size. Therefore, an AP scheme is required. The numerical results for $E$ is showed  in Figure \ref{fig:exfourE} $(a)$.
The results of our new scheme are presented in Table \ref{tab:4} and Figure \ref{fig:exfour}. Similar to all previous examples, the new scheme has uniform second order convergence in space regardless of the anisotropy strength. The numerical
results for different $\epsilon$ are plotted in Figure \ref{fig:exfour2}. When $\epsilon$ becomes small, the solution tends to constant along the $\mathbf{b}$ field line.
\begin{table}[tbp]
\centering
 \begin{tabular}{|cccccc|}
\hline
$ \epsilon\backslash I\times J$  &$32\times32$  & $64\times64$ & $128\times128$ & $256\times256$& $512\times512$   \\ \hline
10 		&1.2232E-02	&3.4612E-03	&7.9300E-04	&2.0647E-04
 &5.4257E-05

 \\ \hline
 1  &3.4397E-03	&8.2860E-04	&2.5018E-04	&6.6258E-05
	&1.7022E-05
\\ \hline
0.5  	&3.1623E-03	&7.9053E-04	&2.3494E-04	&6.2653E-05
 &1.6474E-05
\\ \hline
$10^{-1}$ 	&2.1192E-03	&5.6832E-04	&1.6677E-04	&4.5073E-05
	&1.1970E-05
 \\ \hline
$10^{-3}$  	&1.6690E-03	&4.3141E-04	&1.0753E-04	&2.7379E-05
	&6.9593E-06
\\ \hline
$10^{-6}$ 	&1.6753E-03	&4.3215E-04	&1.0744E-04  &2.7238E-05
	&6.8520E-06
\\ \hline
$10^{-9}$  	&1.6753E-03	&4.3215E-04	&1.0744E-04 &2.7233E-05
	&6.8454E-06
\\ \hline
$10^{-12}$  &1.6746E-03	&4.3311E-04	&1.0795E-04	&2.7403E-05
	&6.9244E-06
\\ \hline
\end{tabular}
 \caption{Example 4. The discrete $L^2$ norm of the errors for different grids and $\epsilon$ values.}
\label{tab:4}
\end{table}

\begin{figure}[htb]
\centering
{
\subfigure [] {\includegraphics[width=7.5cm,clip]{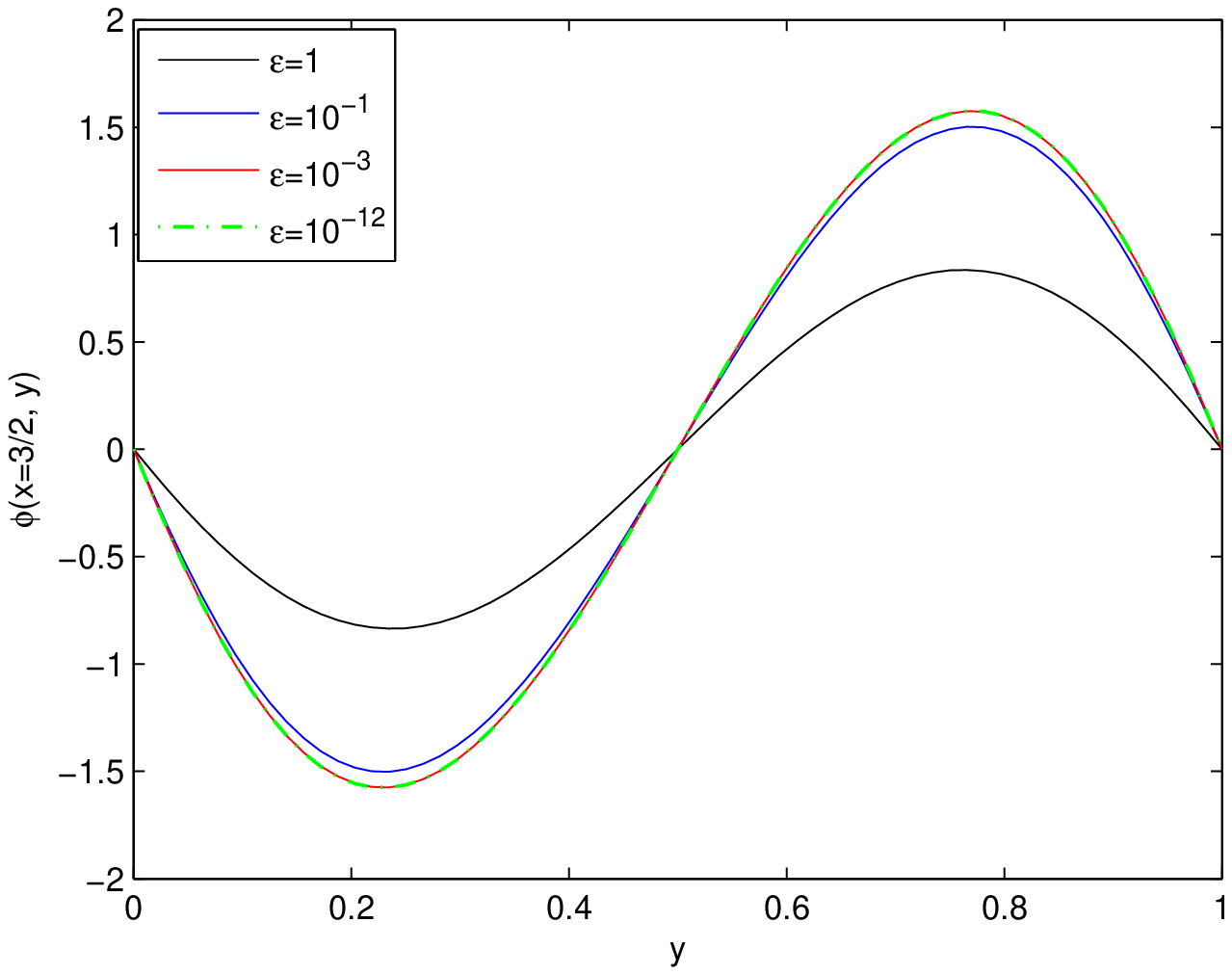} }\\
\subfigure [] {\includegraphics[width=7.5cm,clip]{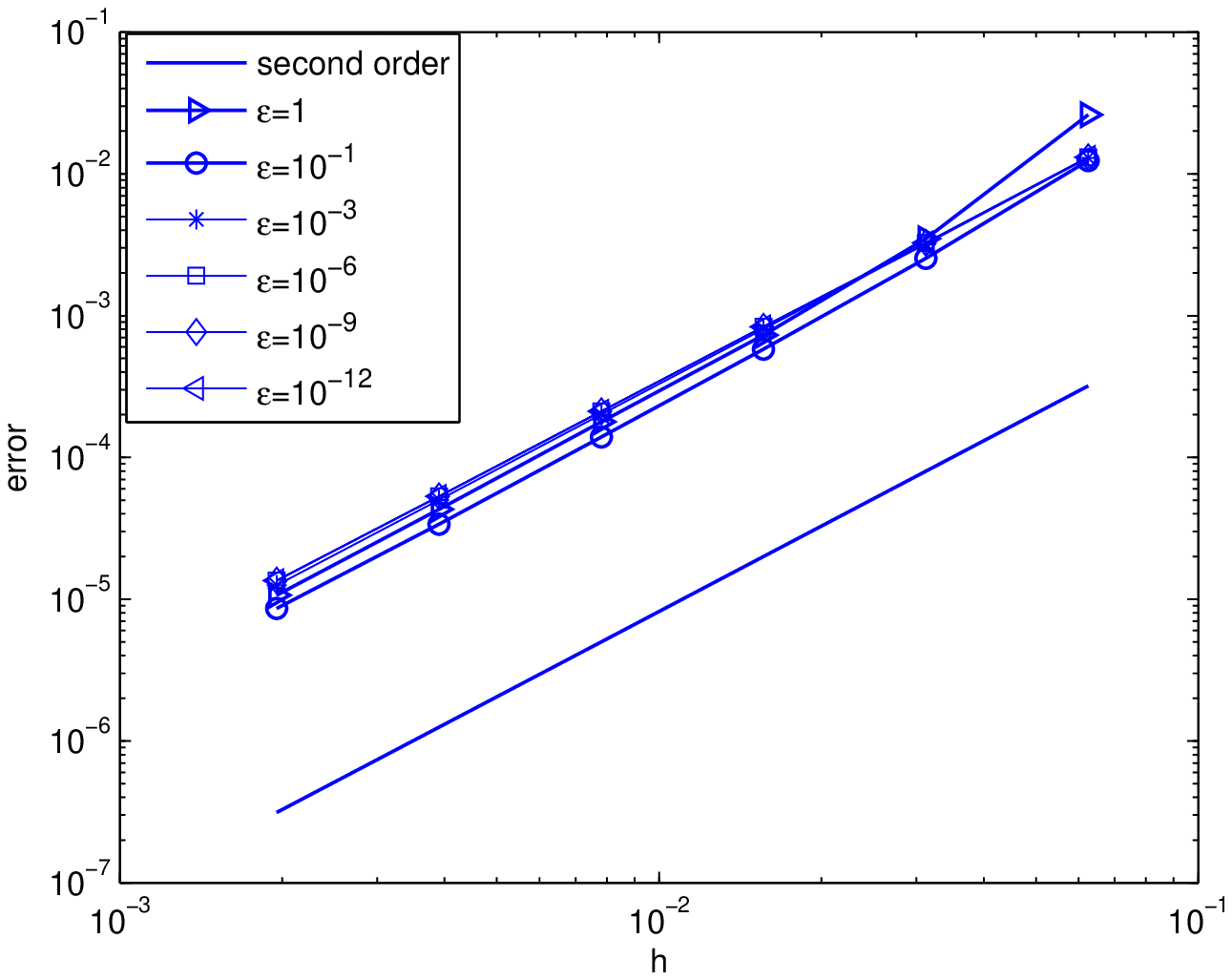} }
\subfigure[] {\includegraphics[width=7.5cm,clip]{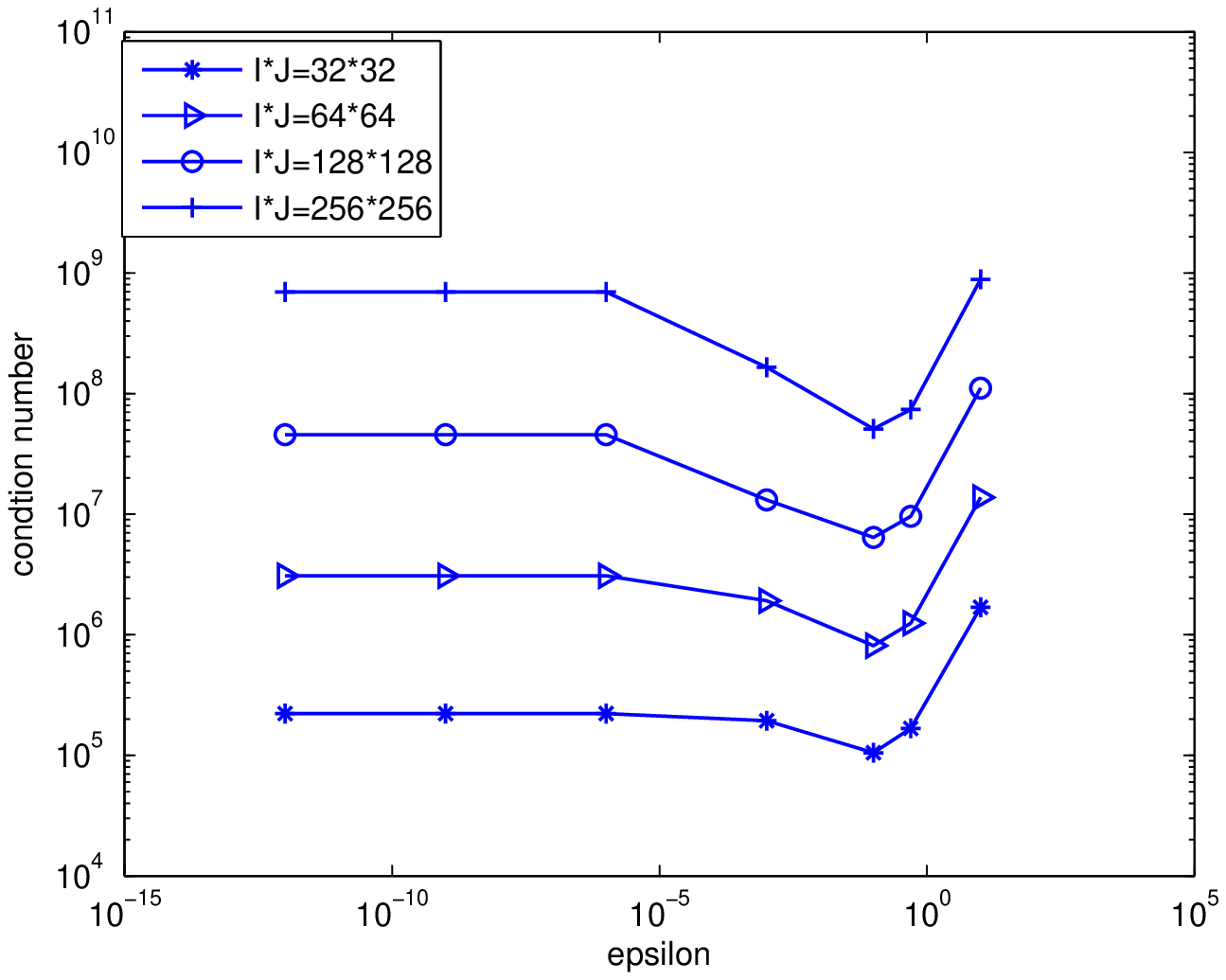}}

}
\caption{Example 4.
The results of the computation domain: $[0,3/2]\times [0,1]$ with different grids for different $\epsilon$ are given. (a) The value of $\phi(x,y)$ in \eqref{eq:test4b} at the right boundary; (b) The convergence order; (c) The condition number.}
\label{fig:newex4}
\end{figure}

 \textbf{Test two:}  The exact solution and the field  are the same as test one, but the computational domain is replaced by   $[0,3/2]\times[0,1]$. Then the Neumann boundary condition
 becomes
 \beq\label{eq:test4b}
 \frac{1}{\epsilon} \mathbf{n}\cdot\mathbf{b}^T  (\mathbf{b}\cdot \nabla u^{\epsilon}) + \alpha \mathbf{n}\cdot\mathbf{b}^T_{\perp}
   (\mathbf{b}_{\perp}\cdot \nabla u^{\epsilon}) = \phi(x,y),
 \eeq
where
$\phi(x,y)$ can be calculated by plugging the exact solution \eqref{eq:test3} into \eqref{eq:test4b}. Thanks to that the field $\mathbf{b}$ satisfies $\mathbf{b}\cdot\nabla u^0=0$,
$\phi(x,y)$ is a bounded function independent of $\epsilon$, as plotted in Figure \ref{fig:newex4} $(a)$.
When $\epsilon\to 0$, the original system \eqref{eq:SimellipA} coupled with the boundary condition \eqref{eq:test4b} at $\Gamma_{out}$ yields the same ill-posed problem as in \eqref{eq:rew2lim}. We replace the boundary condition on $\Gamma_{out}$ by
  the integration of the problem along the field line,
  that is
\beq\label{eq:test2in}
\left[ -E \frac{\phi(x,y)}{\mathbf n \cdot  \mathbf b^T}+E \alpha\frac{ \mathbf n \cdot \mathbf b^T_{\perp}}{\mathbf n \cdot  \mathbf b^T} \mathbf b_{\perp}^{} \cdot
   \nabla u^{\epsilon} \right]_{s = 0}^{s = L_l} - \int_0^{L_l} E \nabla \cdot ( \alpha
   \mathbf b_{\perp}  (\mathbf b_{\perp}\cdot \nabla u^{\epsilon})) {ds} = \int_0^{L_l} Ef{ds}.
\eeq
Replacing the right end of the equation \eqref{eq:discL} by
$$
\Big( E \frac{-\phi(x,y)}{\mathbf n \cdot  \mathbf b^T}+ E\alpha\f{\mathbf{n}\cdot\mathbf{b}^T_\bot}
    {\mathbf{n}\cdot\mathbf{b}^T}\big(-\sin\theta\p_xu^\epsilon+\cos\theta\p_yu^\epsilon\big)\Big)\Big|_{(x_{0},\bar{y}^k_{0})}^{(x_{I},\bar{y}^k_{I})}
$$
 we can get the discretization for \eqref{eq:test2in}.
 The numerical results for $E$ is showed in Figure \ref{fig:exfourE} $(b)$.
The convergence orders and condition numbers for different $\epsilon$  are plotted in Figure \ref{fig:newex4}. The scheme has uniform second order convergence when $\mathbf n \cdot
  \mathbf{b}^T_{\perp} \neq 0$ at $\Gamma_{out}$, regardless of the anisotropy.

\begin{figure}[htb]
\centering
{
\subfigure [] {\includegraphics[width=7.5cm,clip]{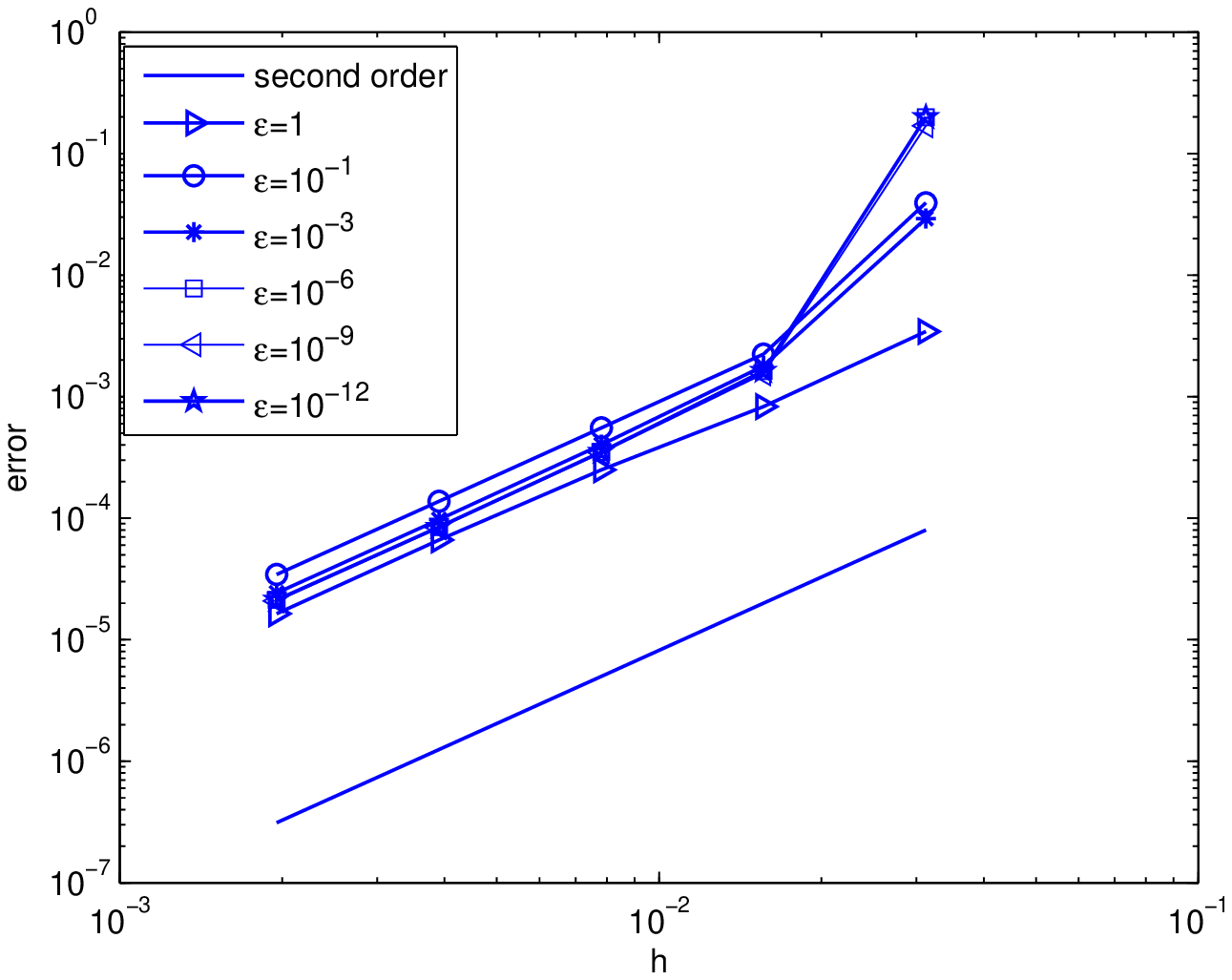} }
\subfigure [] {\includegraphics[width=7.5cm,clip]{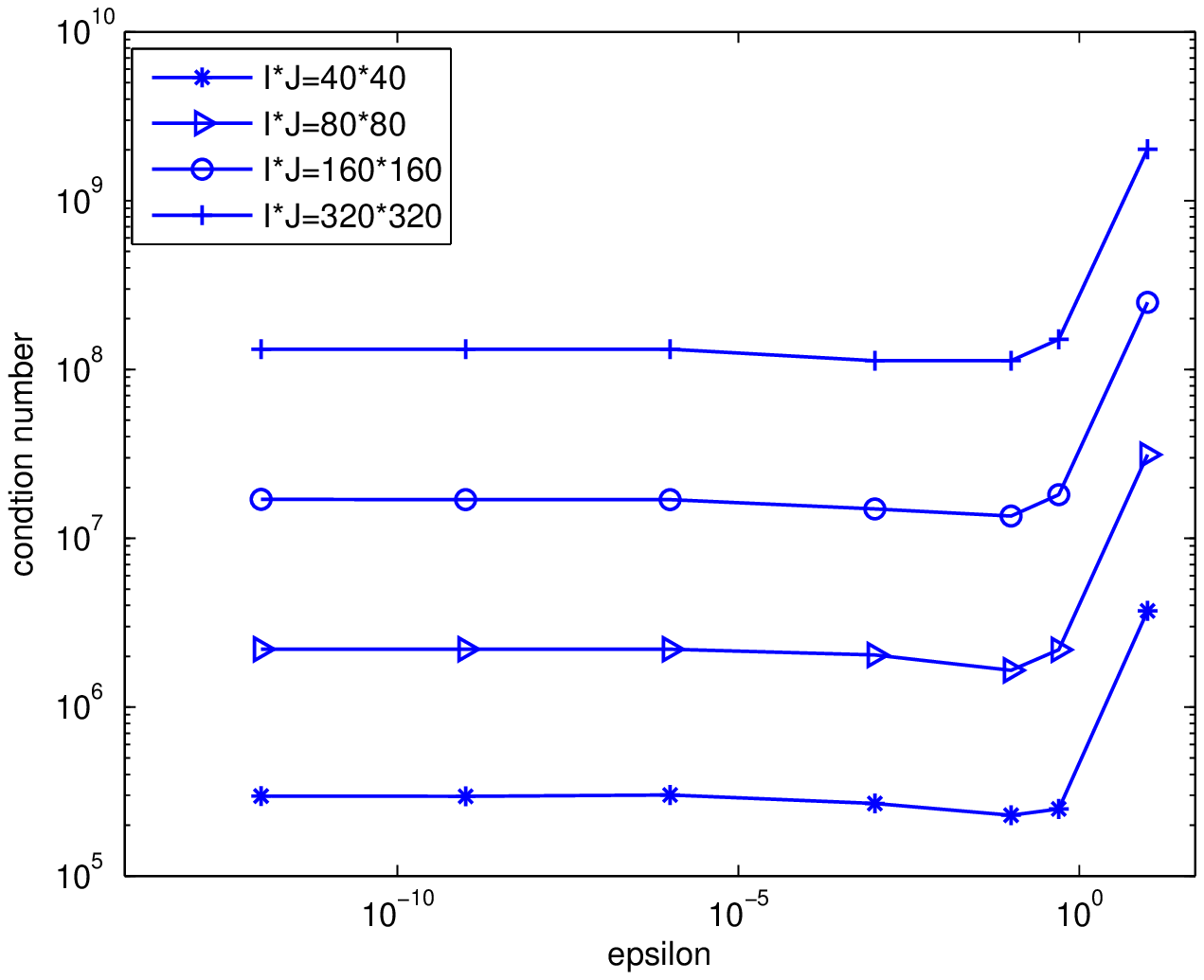} }

}
\caption{Example 5.
  The results with different grids for different $\epsilon$. (a): Convergence orders; (b): Numerical condition number.}\label{fig:exfive}
\end{figure}

\begin{figure}[htb]
\centering
{
\subfigure [] {\includegraphics[width=7.5cm,clip]{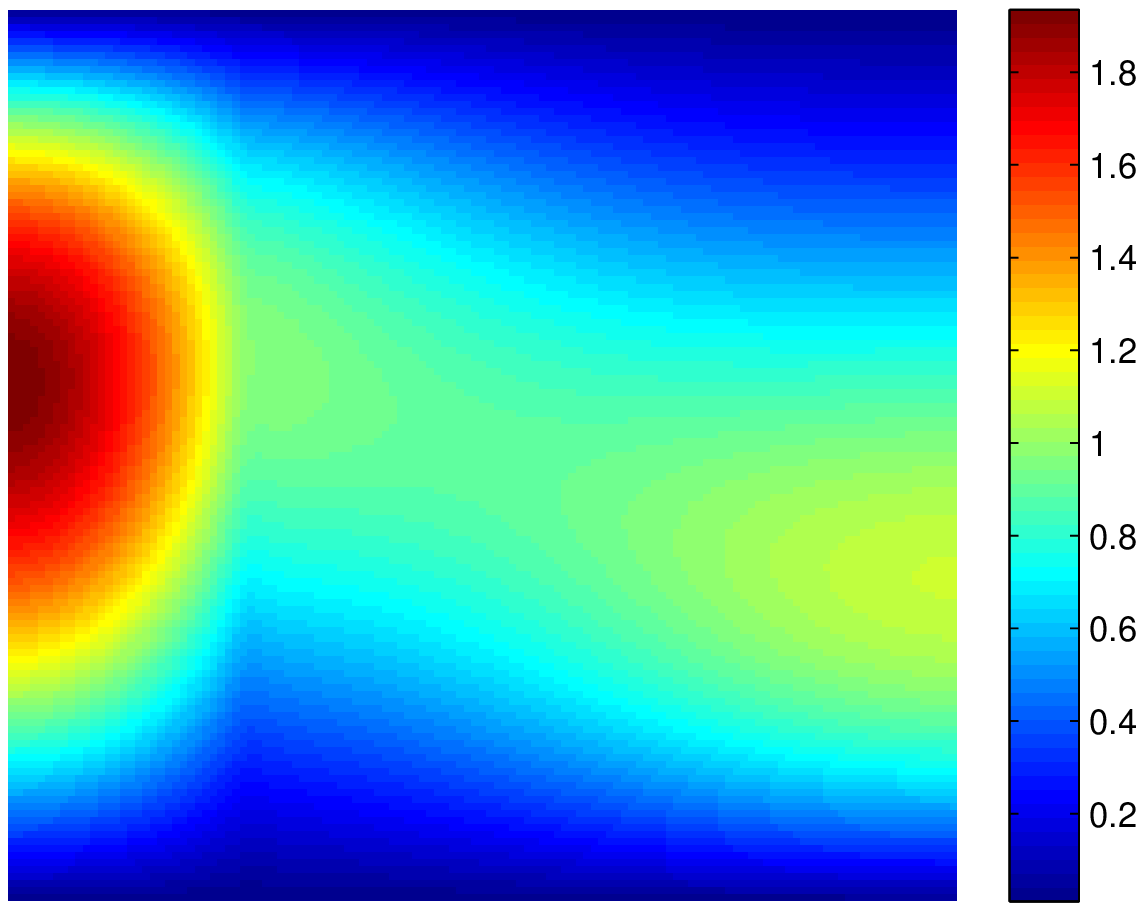} }
\subfigure [] {\includegraphics[width=7.5cm,clip]{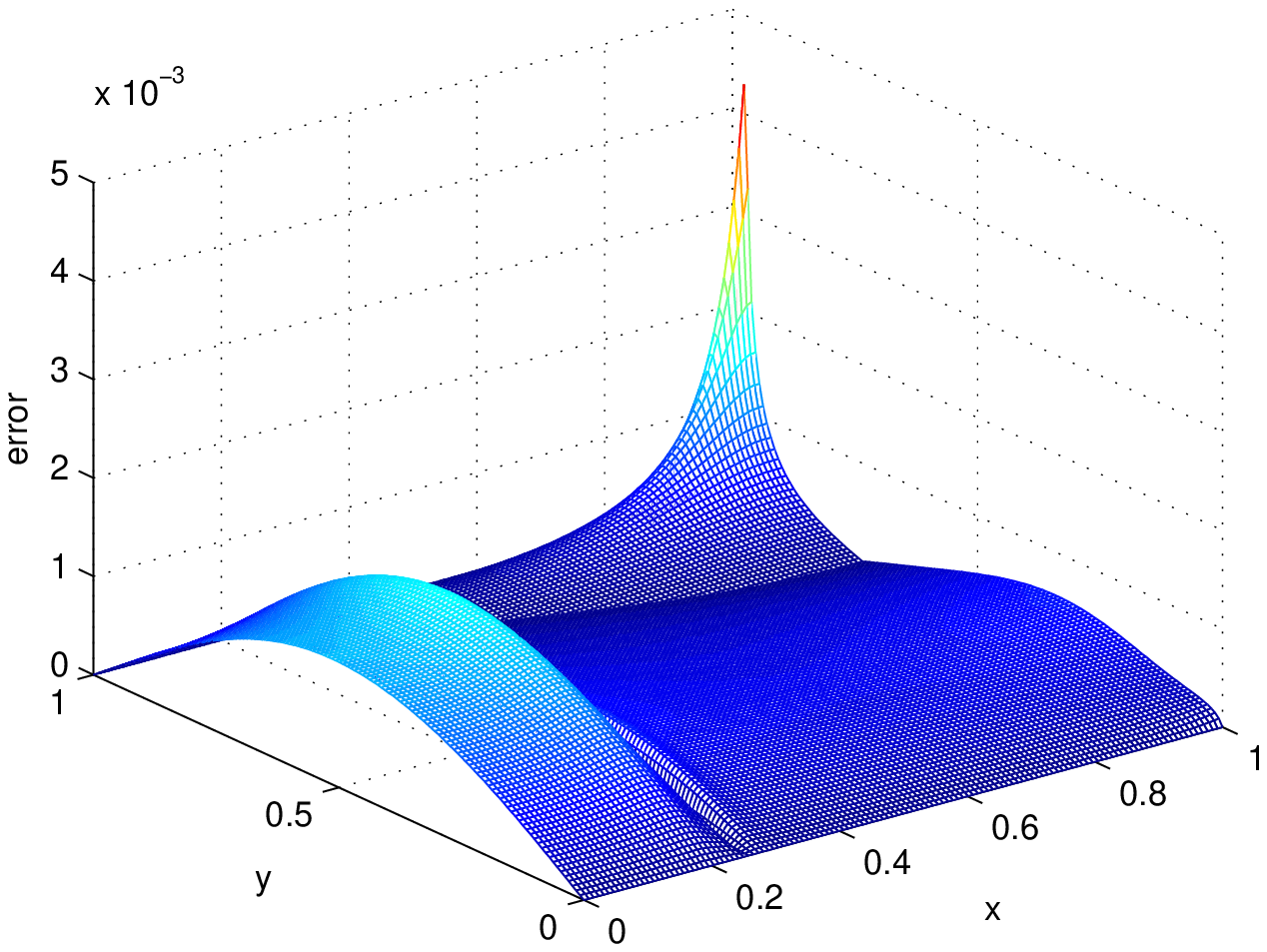} }
\subfigure[] {\includegraphics[width=7.5cm,clip]{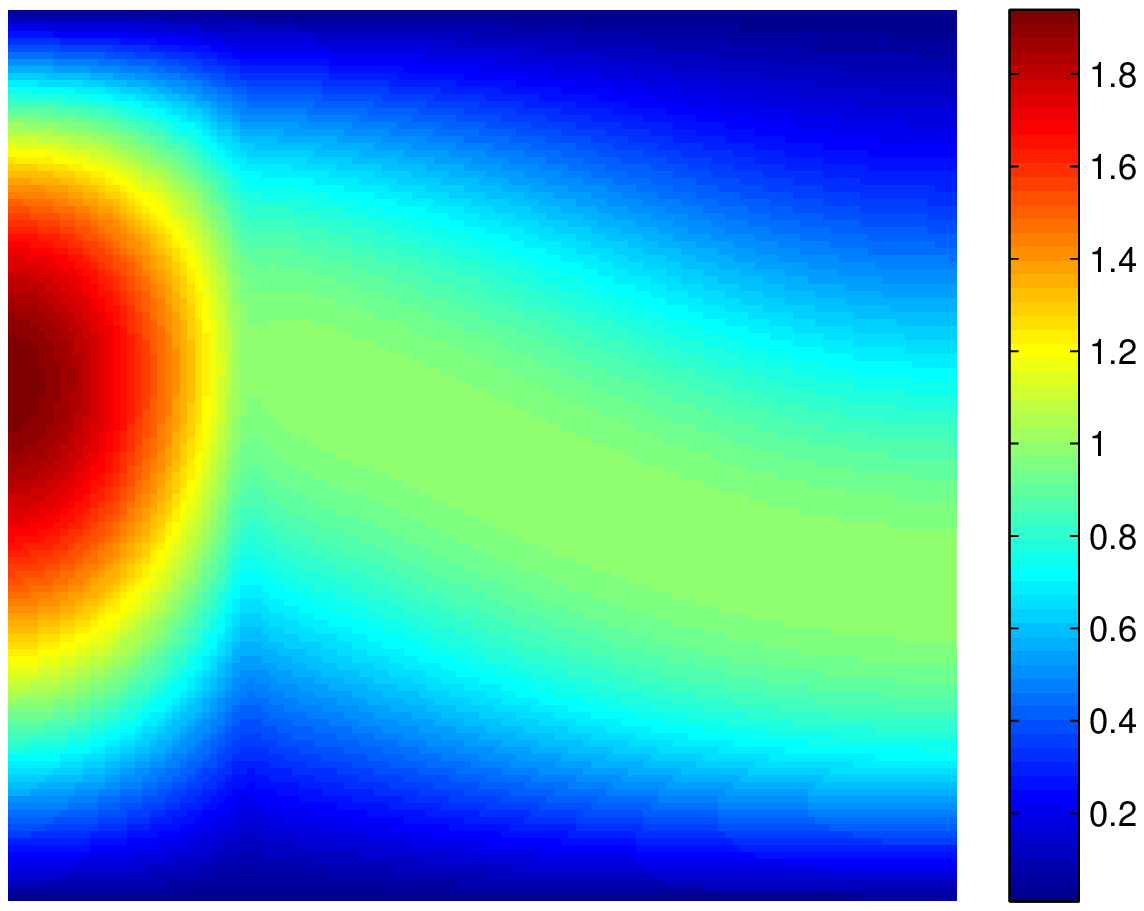}  }
\subfigure [] {\includegraphics[width=7.5cm,clip]{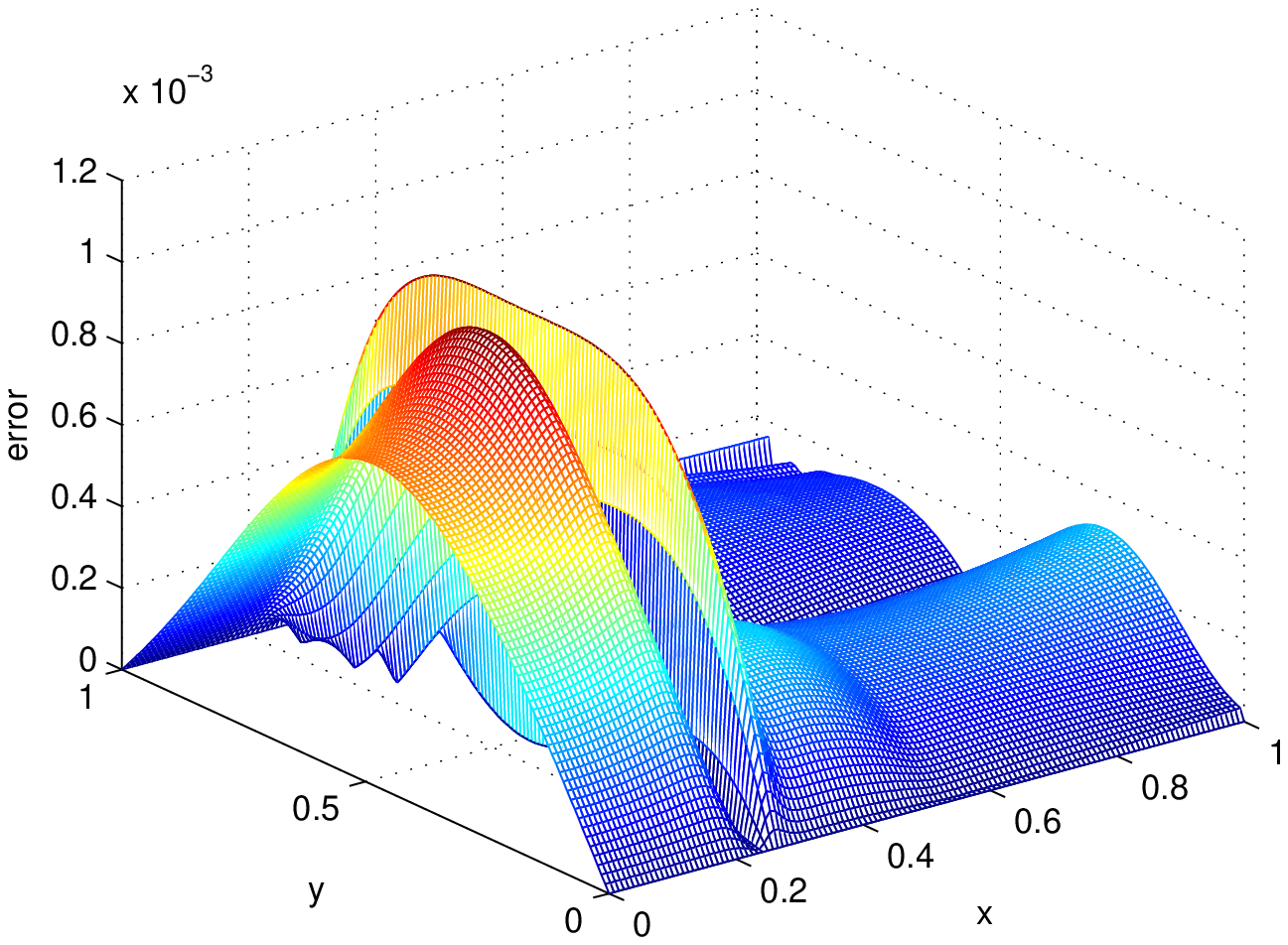} }

}
\caption{Example 5.
 (a): The numerical results of $u^\epsilon(x,y)$ with $\epsilon_{min}=0.1$ in $\epsilon(x,y)$, $I\times J=128\times128$.
 (b): The numerical error with $\epsilon_{min}=0.1$ and $I\times J=128\times128$.
 (c): The numerical results of $u^\epsilon(x,y)$ with $\epsilon_{min}=10^{-12}$ and $I\times J=128\times 128$.
 (d): The numerical error with $\epsilon_{min}=10^{-12}$ and $I\times J=128\times128$.}\label{fig:exfive2}
\end{figure}

\bigskip
\textbf{Example 5:} Variable $\epsilon$ and non-aligned $\mathbf{b}$ field case. We choose the same limiting solution as well as the magnetic fields as in Example 4,
 but $\epsilon$ varies in the computational domain and is defined by \eqref{eq:vepsilon}. The analytical solution of the present example is
 \beq\label{eq:test5}
u_{exact}(x,y)=\sin(\pi y+\kappa(y^2-y)\cos(\pi x))+\epsilon(x,y)\cos(2\pi x)\sin(\pi y),
\eeq
which satisfies the Neumann boundary condition on the left and right boundary of the computational domain $[0,1]\times[0,1]$.
 The source term is calculated accordingly.  The performance of our new scheme is presented in Table \ref{tab:5} and Figure \ref{fig:exfive}.
 The new scheme is capable of producing accurate results when the strength of anisotropy varies
 a lot in the computational domain, even for non-aligned  $\mathbf{b}$ field. As can be seen in Figure \ref{fig:exfive2},  when $\epsilon_{min}$ become small, $u^\epsilon$ is almost constant along the field lines
 at the part $\epsilon(x,y)$ is small, while the largest error appear at the
 transition region. {Similar as the discussions in Example 3, since the mesh size of $I \times J = 32 \times 32$ is too large to accurately represent the narrow transition region, which, therefore, introduces the fast error decreases in Table 5 between the columns $I \times J = 32 \times 32$ and $64\times 64$ for $\epsilon \leq 10^{-1}$.}

\begin{table}[tbp]
\centering
\begin{tabular}{|cccccc|}
\hline
$ \epsilon\backslash I\times J$     &$32\times32$  & $64\times64$ & $128\times128$ & $256\times256$& $512\times512$   \\ \hline
10 	&8.7766E-02		&8.9015E-03		&2.2852E-03		&5.7277E-04		&1.4330E-04
 \\ \hline
1  		&3.4397E-03	&8.2860E-04	&2.5018E-04&6.6258E-05
	&1.7022E-05
\\ \hline
0.5 	&6.6596E-02		&8.3994E-04		&2.4285E-04		&6.4713E-05  &1.6182E-05
\\ \hline
$10^{-1}$  &3.9386E-02  &2.2373E-03	&5.5221E-04		&1.3778E-04   & 3.4361E-05
 \\ \hline
$10^{-3}$  	&2.9133E-02		&1.7894E-03		&4.0130E-04		&9.7669E-05 &2.4210E-05
\\ \hline
$10^{-6}$  	&1.9871E-01		&1.6253E-03		&3.5132E-04		&8.4812E-05 &2.1265E-05
\\ \hline
$10^{-9}$  	&1.6929E-01		&1.5458E-03		&3.4572E-04		&8.3489E-05  &2.0843E-05
\\ \hline
$10^{-12}$	&1.2154E-01		&1.5791E-03		&3.4117E-04		&8.4430E-05  &2.1271E-05
\\ \hline
\end{tabular}
 \caption{Example 5. The discrete $L^2$ norm of the errors for different grids and $\epsilon$ values.}
\label{tab:5}
\end{table}

 \section{Conclusion}\label{con}
We present a simple Asymptotic-Preserving reformulation for strongly anisotropic problems.
The key idea is to replace the Neumann boundary condition on one side of the field line by the integration of the original problem along the
 field line. The new system can remove the ill-posedness in the limiting problem and the reason is illustrated at both continuous and discrete level.
 First order methods for the field line integration is used to illustrate the idea but the convergence order can be easily improved by interpolation and trapezoidal
  rule.  Uniform second order convergence can be observed numerically and the condition number does not scaled with
strength of the anisotropy.

The scheme is efficient, general and easy to implement.  Small modifications to the original code are required, which makes it attractable to engineers.
The idea can be coupled with most standard discretizations and the computational cost keeps almost the same.
We can further improve the convergence order by higher order space discretization and curve integration.

 The AP  Macro-Micro decomposition methods introduced in \cite{Degond101,Degond102,Degond121,Degond122} have some difficulties in case of closed field lines. Narski and Ottaviani \cite{Narski13} introduced a penalty stabilization term to improve the accuracy in case of closed field lines, where the penalty stabilization has a tuning parameter and a system
 of two equations instead of one is solved. Extensions of our idea to closed field lines seem straight forward.
Its applicability to closed field lines, time dependent problems as well as nonlinear problems (as in \cite{Mentrelli12}) will be our future work.

\section*{Acknowledgement}
The authors are partially supported by NSFC 11301336
and 91330203.

\end{document}